\documentstyle[12pt,twoside, amscd]{amsart}

\newcommand{\cM}{{\cal M}}

\def\mapdown#1{\Big\downarrow\rlap{$\vcenter{\hbox{$\scriptstyle#1$}}$}}
\def\mapright#1{\smash{   \mathop{\longrightarrow }\limits^{#1}}}
\def\mapleft#1{\smash{   \mathop{\longleftarrow }\limits^{#1}}}
\def\mapdown#1{\Big\downarrow
    \rlap{$\vcenter{\hbox{$\scriptstyle#1$}}$}}
\def\mapup#1{\Big\uparrow
    \rlap{$\vcenter{\hbox{$\scriptstyle#1$}}$}}
\def\proof{\smallskip\noindent{\it Proof. }}
\def\endproof{\hfill\qed}
\def\qed {\nobreak$\quad$\lower 1pt\vbox{
    \hrule
    \hbox to 8pt{\vrule height 8pt\hfil\vrule height 8pt}
      \hrule}\ifmmode\relax\else\par\medbreak \fi}
\def\exact#1#2#3{0\longrightarrow{#1}\longrightarrow{#2}\longrightarrow{#3}\longrightarrow0}

\newtheorem{thm}{Theorem}[section]

\newtheorem{lem}[thm]{Lemma}
\newtheorem{cor}[thm]{Corollary}

\newtheorem{prop}[thm]{Proposition}

\theoremstyle{definition}
\newtheorem{defn}[thm]{Definition}
\newtheorem{condition}[thm]{Condition}
\newtheorem{say}[thm]{}
\newtheorem{exmp}[thm]{Example}

\newtheorem{rem}[thm]{Remark}           %\renewcommand{\theremark}{}
\newtheorem{notationnum}[thm]{Notation}
\theoremstyle{remark}

\begin{document}
\title[Birational Models of the Moduli of Stable Bundles]
{Birational Models of the Moduli Spaces of Stable Vector Bundles over Curves}
\author[Yi {\sc Hu} and Wei-Ping Li]{Yi {\sc Hu}$^*$ and Wei-Ping Li$^{**}$}
\date{November, 1996} %(Preliminary Version)}
%\subjclass
%Primary 14D25, 14L30,  Secondary 14C20, 14A99 \endsubjclass
\thanks{$^*$Research partially supported by NSF grant DMS 9401695.} 
\thanks{$^{**}$Research partially supported by RGC Earmarked Grant 
HKUST631/95P and DAG Grant from HK government.}
\begin{abstract}
We give a method to construct stable vector bundles
whose rank divides the degree over curves of genus bigger than one.
The method complements the one given by Newstead.
Finally, we make some systematic remarks and observations 
in connection with rationality of moduli spaces of stable vector bundles.
%We also show that Newstead's induction method can be used to prove
%rationality of a larger class of moduli spaces of stable vector bundles
%over curves of genus bigger than one. 
\end{abstract}

\maketitle

{\footnotesize
%{\scriptsize
\tableofcontents
}

\section{Introduction}
\label{sec:introduction}
\begin{say} 
Let $X$ be a complete non-singular 
algebraic curve of genus $g \ge 2$ over an algebraically closed field
$k$ of characteristic 0 (for simplicity we assume that $k$ is the field
of complex numbers ${\Bbb C}$). The moduli space ${\cal M}(n, L)$ of semi-stable
 vector bundles of rank $n$ and determinant 
$L$ with degree ${\rm deg} L=d$ over $X$ 
is an irreducible projective variety of 
dimension $(n^2 -1)(g-1)$. The question whether the moduli space ${\cal M}(n, L)$ is rational
is  very subtle. 
The affirmative answer is known for many cases
through the works of 
 Narasimhan-Ramanan ([{\bf NR}]), Newstead ([{\bf N1, N2}]), 
and Tjurin ([{\bf T3}]). 
In proving rationality in [{\bf N1, N2}], 
Newstead  conjectured  a systematic way to construct 
(generic) stable vector bundles
$F$ from (generic) stable vector bundles $F'$ of lower rank via extensions of
the following type:
$$\exact{{\cal O}_X^{\oplus r}}{F}{F'}$$
where ${\rm rank}(F)=n$ and ${\rm deg}F=d=n(g-1)+r$ for some $0<r<n$. 
He proved the conjecture in [{\bf N1, N2}] for some cases and
Grzegorczyk completed the proof for all cases in her paper [{\bf G}].
In [{\bf N2}], this extension was used as an induction
step to show that the following moduli spaces of stable vector bundles
are rational.
\end{say}

\begin{thm} ({\rm [{\bf N2}]})
\label{newsteadrationalitythm}
The moduli space ${\cal M}(n, L)$ is rational in the following
cases:
\begin{enumerate}
\item  $d=1$ mod $n$ or $d=-1$ mod $n$;
\item $\gcd (n, d)=1$ and $g$ is a prime power;
\item $\gcd(n, d)=1$ and the sum of the two smallest distinct prime factors of $g$ is greater than $n$.
\end{enumerate}
\end{thm}

\begin{say} Clearly, from that $d=n(g-1)+r$ and $0<r<n$, 
one sees that the above method of constructing stable vector bundles leaves out the
case when rank divides degree. {\sl It is one of the objectives of the current paper to find a method of 
constructing stable vector bundles in complement to that of} [{\bf N1, N2, G}].
\end{say}

\begin{thm}
\label{constructingstablebundles}
Let $n$, $d$ be positive integers such that 
$d=ng$. Let $L={\cal O}_X(P_0)$ where $P_0$ is a special effective divisor
defined in \ref{choosedivisor}.
Then there exists a non-empty Zariski open subset of the moduli
space ${\cal M}(n, L)$ consisting of vector bundles $V$ such that
\begin{enumerate}
\item $h^1(X, V)=0$.
\item ${\cal O}_X^{\oplus n}$ is a sub-sheaf of $V$.
\item there exists an exact sequence 
$$0\longrightarrow V^*\longrightarrow {\cal O}_X^{\oplus n}\mapright{\varphi}
 {\cal O}_P\longrightarrow 0\eqno(1)$$
where $P$ is a divisor $P=p_1+\ldots +p_d$ such that the map
 ${\varphi}$
in $Hom({\cal O}_X^{\oplus n}, {\cal O}_P)$ 
is stable with respect to the action of $Aut({\cal O}_X^{\oplus n}) \times Aut({\cal O}_P)$.
\end{enumerate}
\end{thm}

\begin{rem} We point out that Tjurin ([{\bf T1, T2}]) studied these moduli spaces using 
the different method of matrix divisors.
\end{rem}

For the purpose of comparison, we mention a  result  of 
Grzegorczyk [{\bf G}] who proved the following conjecture of Newstead:

\begin{thm} {\rm  ([{\bf N1, G}]) }
\label{gthm}
Let, $n$, $d$, $r$, be natural numbers such that 
$d=n(g-1)+r$, $0<r<n$. Then there exists a non-empty Zariski open subset
of the moduli space ${\cal M}(n, L)$ consisting of vector bundles $F$ such that 
\begin{enumerate}
\item $h^1(F)=0$.
\item ${\cal O}_X^{\oplus r}$ is a sub-bundle of $F$.
\item the quotient bundle $F/{\cal O}_X^{\oplus }$ is stable.
\end{enumerate}
\end{thm}

This theorem covers all moduli spaces except those where $n\, |\, d$. Our
theorem above deals exactly the complement.

\begin{say}
The tool that we used to construct stable vector bundles
is that of elementary transformations which were probably
first studied by Maruyama ([{\bf M}]). We choose
a special effective divisor $P_0$ (see \ref{choosedivisor})
whose degree equals $ng$. Take $L$ to be
${\cal O}_X(P_0)$. Let $U$ be a non-empty Zariski open subset of the linear system
$|P_0|$ satisfying certain generic condition (see \ref{chooseU}). Let $P$ be an 
effective divisor in $U$. Now given a surjective map $\varphi\in
Hom({\cal O}_X^{\oplus n}, {\cal O}_P)$ (which is called an elementary transformation),
we get an exact sequence
$$0\longrightarrow {W}\longrightarrow {{\cal O}^{\oplus n}_X}\mapright{\varphi}
{\cal O}_P\longrightarrow 0.$$
Let $V=W^*$. Then we see that ${\rm rank} V=n,\quad {\rm det} V={\cal O}_X(P)={\cal O}_X(P_0)=L.$
Then Theorem \ref{constructingstablebundles} (see also \S\S 3--5) basically
 asserts that if the elementary transformation 
satisfies some generic conditions, then
$V$ is stable and $h^1(V)=0$. Furthermore, the
stable vector bundles $V$ obtained this way form a non-empty
Zariski open subset of the moduli space and hence it leads to 
a birational model for the moduli space.
\end{say}

\begin{say} In the end of this paper, we shall systematically explore that 
the inductive method in ([{\bf N2}]) can actually be extended to
prove rationality of a rather larger class of moduli spaces.
This class consists of  much more cases than those listed
 in Theorem \ref{newsteadrationalitythm}.

For example, when $g=6$, $n=15, d=77$, the moduli space ${\cal M}(n, L)$ is rational.
Also, the moduli spaces ${\cal M}(n, L)$ with 
$$n =11+7m, \quad  {\rm deg}L=62+35m$$ 
are  all rational over any genus $6$ curve
for $m\ge 7$. But the set $\{ g=6, n=15, d=77 \}$ and
$\{ g=6, n=11+7m, d= 62+35m \}$ satisfy none of the three
conditions as listed in Theorem \ref{newsteadrationalitythm}. 

In general, the moduli space ${\cal M}(n, L)$ is 
rational if $n(g-1) < d={\rm deg}L < ng$ and the pair $(n; d)$
can be ``linked'' to $(1; g)$ by some successive arithmetic reductions of
the following types:
\begin{enumerate}
\item  $(n'; d')=(ng-d; d-k(ng-d))$ for some non-negative integer $k$ such that
       $n'(g-1) < d' \le  n'g$; or
\item $(n''; d'')=(d-n(g-1); n(2g-1)-d-k(d-n(g-1)))$
 for some non-negative integer $k$ such that
       $n''(g-1) < d'' \le  n''g$. 	
\end{enumerate}
This fact (Theorems \ref{newsteadpairthm} and \ref{ballicopairthem}), 
 together with some other results, is systematically explored in \S 6.
%was basically pointed out by Newstead  using properties of the 
%corresponding  moduli spaces. But the above involves only arithmetics and 
%should be a better interpretation of his result. 

Also, we shall point out how our Theorem \ref{constructingstablebundles} is related to
other moduli spaces by similar arithmetic reductions as above.
\end{say}

\begin{say} The paper is structured as follows.
 In \S 2, we recall the Abel-Jacobi theory and other technical results that will
be used in the latter sections.
\S 3 deals with elementary 
transformations and how to use them to  construct stable vector bundles as stated
in Theorem \ref{constructingstablebundles}. \S 4 studies the group action
of $Aut({\cal O}_P)$ and $Aut({\cal O}_X^{\oplus n})$ on $Hom({\cal O}_X^{
\oplus n}, {\cal O}_P)$. \S 5 concerns
birational models for the moduli spaces and completes the proof of our main result
Theorem \ref{constructingstablebundles}. 
Finally, in \S 6,  we make some systematic remarks and observations on rationality of
moduli spaces.
\end{say}

\begin{say} We fix the following notations:
\par\noindent 
$V$ is a locally free sheaf, by abuse of notation, we also call it
a vector bundle.
\par\noindent $V^*$ is the  dual of $V$, i.e. 
$V^*={\cal H}om({\cal O}_X, V)$.
\par\noindent $H^i(V)$ is the cohomology $H^i(X, V)$.
\par\noindent $h^i(V)$ is the dimension of $H^i(V)$.
\par\noindent
$n$ is the rank of the vector bundle $V$.
\par\noindent
$d$ is the degree of the vector bundle $V$.
\par\noindent
${\rm det} V=\wedge^n V$ is the determinant line bundle.
\par\noindent
$L$ is a line bundle with  ${\rm deg} L=d$.
\par\noindent $(M)^d= M\times\ldots\times M$ ($d$-times)
%$(M)^d=\underbrace{M\times\ldots\times M}_d.$
\par\noindent (2.1), (2.2), ... ,  are the indices labeling formulas.
\par\noindent 2.1, 2.2, ... ,  are the indices labeling theorems, remarks, paragraphs, and so on.
\end{say}

This note stems from our study of the rationality problem 
of moduli spaces  when both of the authors were visiting
Max-Planck-Institute f\"ur Mathematik in the summer of 1993.
Due to the fact that there are few papers in the literature addressing the moduli spaces when
the rank divides the degree, we decided to write up this note and make it
available to interested readers.

\medskip\noindent
{\bf Acknowledgments:}
The hospitality and financial support from the 
Max-Planck-Institute (summer, 1993) are  gratefully acknowledged. 
We thank P.  Newstead, Igor Dolgachev, A. Beauville for their attention on this work.
 W.L. thanks  M. Reid and MRC of University of Warwick for
providing a stimulating environment for research. 
Y.H. acknowledges the Centennial Fellowship award 
by the American Mathematical Society.

\section{Preliminaries}

\begin{defn}
Let $X$ be a smooth projective curve over the field of complex numbers
and $V$ a rank-$n$ algebraic vector bundle over $X$.
$V$ is stable (semi-stable) if for any proper sub-vector
 bundle $F$ of $V$, 
$${{\rm deg }F\over {\rm rank }F}<\,\,(\le)\,\,{{\rm deg}V\over {\rm rank V}}.
\eqno(2.1)$$
 Since any torsion free coherent sheaf over a curve is locally free,
we can get a slightly 
modified definition of stability ( or semi-stability): $V$ is stable ( or semi-stable)
if for any proper sub-vector bundle $F$ of $V$ whose cokernel $V/F$ is a vector
bundle, (2.1) holds.
\end{defn}

\begin{say} 
Let $L$ be  a line bundle  over $X$. Throughout this paper,
${\cal M}(n, \,L)$ represents the moduli space of rank-$n$ semi-stable 
vector bundles over $X$ with ${\rm det}V
\buildrel\sim\over =L$.
\end{say}

\begin{say} 
\label{choosedivisor}
We begin with  a special 
divisor whose corresponding line bundle
will be chosen as the fixed determinant of
our semi-stable bundles.

Let $\omega_1,\ldots, \omega_g$ be a basis of $H^0(K_X)$. By the Abel-Jacobi theory
(see [{\bf GH}]), for general points $x_1,\ldots, x_g\in X$, the determinant
$$ \left| \begin{array}{clcr}
\omega_1(x_1) & \omega_1(x_2) &  \ldots & \omega_1(x_g) \\
\vdots & \vdots & \ddots &  \vdots  \\
\omega_g(x_1)&  \omega_g(x_2) & \ldots & \omega_g(x_g) 
\end{array} 
\right| $$ 
is not zero. 
 We now choose, once and for all,  $d=ng$ many distinct
 points $q_1, \ldots, q_{ng}$ such that the determinants
$$\left| \begin{array}{cccc}
\omega_1(q_{i+1})&\omega_1(q_{i+2})& \ldots & \omega_1(q_{i+g}) \\
\vdots&\vdots&\ddots&\vdots \\
\omega_g(q_{i+1})&\omega_g(q_{i+2})&\ldots&\omega_g(q_{i+g}) 
\end{array} \right| \neq 0 \eqno(2.2) $$
for  $i=0, g, 2g,  \ldots, (n-1)g$. Here 
 $n$ is a positive integer bigger than one which will be taken 
later on as the rank of vector bundles. Let $P_0$ be the effective divisor 
$P_0=q_1+\ldots+q_{ng}$. The linear
system $|P_0|$ is a projective space of dimension $(n-1)g$. This can be easily seen 
by Riemann-Roch.
\end{say}

\begin{say}
For a technical reason as we shall see in \S 3, we need to consider the determinant
$$\left| \begin{array}{cccc}
b_{11}\omega_1(x_1)&b_{21}\omega_1(x_2)&\cdots&b_{d1}\omega_1(x_d)\\
b_{11}\omega_2(x_1)&b_{21}\omega_2(x_2)&\cdots&b_{d1}\omega_2(x_d) \\ 
\vdots&\vdots&\ddots&\vdots \\ 
b_{11}\omega_g(x_1)&b_{21}\omega_g(x_2)&\cdots&b_{d1}\omega_g(x_d) \\ 
b_{12}\omega_1(x_1)&b_{22}\omega_1(x_2)&\cdots&b_{d2}\omega_1(x_d) \\ 
\vdots&\vdots&\ddots&\vdots \\ 
b_{12}\omega_g(x_1)&b_{22}\omega_g(x_2)&\cdots&b_{d2}\omega_g(x_d) \\ 
\vdots&\vdots&\ddots&\vdots \\ 
b_{1n}\omega_1(x_1)&b_{2n}\omega_1(x_2)&\cdots&b_{dn}\omega_1(x_d) \\ 
\vdots&\vdots&\ddots&\vdots \\ 
b_{1n}\omega_g(x_1)&b_{2n}\omega_g(x_2)&\cdots&b_{dn}\omega_g(x_d) 
\end{array} \right|. \eqno (2.3)$$
Take $b_{ij}$'s as unknowns, the determinant at the point 
$(q_1,\cdots, q_d)\in (X)^d$ (chosen as in \ref{choosedivisor})
is a non-vanishing polynomial in variables $b_{ij}$. The way to see this is as follows.
We expand the determinant of (2.3) and check the coefficient of 
$$b_{11}b_{21}\ldots b_{g1}b_{(g+1)2}\ldots b_{(2g)2}\ldots b_{(d-g)n}\ldots
b_{dn}.$$
By the technique of 
minor expansions
 in the determinant theory,  the coefficient is
$$\hbox{$\prod\limits_{i=0}^{(n-1)g}$}
\left| \begin{array}{cccc} 
\omega_1(q_{i+1})&\omega_1(q_{i+2})& \ldots&\omega_1(q_{i+g}) \\ 
\vdots&\vdots&\ddots&\vdots \\ 
\omega_g(q_{i+1})&\omega_g(q_{i+2})&\ldots&\omega_g(q_{i+g}) 
\end{array} \right| $$
which is not zero by the choices of $q_i$'s.
Hence the determinant 
(2.3) as a polynomial in variables $b_{ij}$ is not identically zero.
\end{say}

\begin{say}
\label{chooseU}
Choose a   collection of $b_{ij}$'s such that the determinant at 
$(q_1,\cdots, q_d)$ is not zero. With these $b_{ij}$'s fixed, 
the zero locus of the determinant (2.3) defines
a divisor in $|P_0|$. %Notice that the value of the determinant (1.3) at the point
%$P\in |P_0|$ is determined up to a sign. 
So we can choose a non-empty Zariski open subset $U$ of $|P_0|$ such that
every divisor $P=p_1+\dots+p_d$ in $U$ satisfies the property that
$p_i$'s are distinct and that the determinant (2.3) with 
the chosen $b_{ij}$'s is non-zero at $P\in U$. 
\end{say}

\begin{say}
It is known that ${\cal M}(n, L) \cong {\cal M}(n, L')$
if ${\rm deg} L={\rm deg} L'$. This can be seen as follows.
There exists a line bundle $\bar L\in Pic^0(X)$ such that
 $\bar L^{\otimes n}=L^*\otimes L'\in Pic^0(X)$. We have that  ${\rm det} (V\otimes
\bar L)={\rm det} V\otimes \bar L^{\otimes n}=L\otimes L^*\otimes L'=L'={\rm det}V'$.
Hence there is a bijection $${\cal M}(n, L)\mapright{\otimes \bar L}
{\cal M}(n,L').\eqno(2.4)$$ 
Since these two moduli spaces are coarse, the map (2.4) must be an
isomorphism. Therefore,  we have some freedom in  choosing a line bundle $L$ in a way
 we like without affecting the isomorphic type of the moduli space ${\cal M}(n, L)$.
In this paper,  we choose $L={\cal O}_X(P_0)$.
\end{say}

\section{Stable vector bundles whose ranks divide their degrees}

This section is entirely devoted to the case when the rank of the bundles divides the degree.
We shall use elementary transformations to construct generic
rank-$n$ stable vector bundles $V$ with degree $ng$ whose 
$H^1(V)$ is trivial. In the section that follows, we shall discuss the necessary
group actions on the space of elementary transformations, which will lead to
a birational model for the moduli space.

\begin{defn} Let $W_1$ $W_2$ be two rank-$n$ vector bundles
over an algebraic variety. Suppose that there exists an injective morphism from
$W_1$ to $W_2$, then we can write a short exact sequence
$$0\longrightarrow{W_1}\longrightarrow {W_2}\mapright{\varphi}{Q}
\longrightarrow 0$$
where $Q$ is a torsion sheaf. In this case, we say that $W_1$ is an elementary 
transformation of $W_2$ with respect to the surjective map $W_2\mapright{\varphi} Q
\longrightarrow 0$; or we may simply say that the map $\varphi$ is an elementary
transformation. %The space $Hom(W_2, Q)$ will be called the space of elementary 
%transformations.  
\end{defn}

\begin{say}  
\label{elemtrans}
Let $P=p_1+\ldots +p_d$  be an effective divisor
over the curve $X$ such that $P\in U$ (see \ref{chooseU}).
Consider an elementary transformation:
$$0\longrightarrow {W}\longrightarrow{{\cal O}_X^{\oplus n}}\mapright{\varphi} {{\cal O}_P}
\longrightarrow 0.  $$
$W$ has to be a vector bundle.  Let $V=W^*$ or
$W=V^*$. Then the above exact sequence can be rewritten as
$$0\longrightarrow {V^*}\longrightarrow
{{\cal O}_X^{\oplus n}}\mapright{\varphi}{{\cal O}_P}\longrightarrow 0.
\eqno(3.1)$$
A simple computation leads to 
$${\rm rank}V=n,\quad {\rm det}V={\cal O}_X(P)={\cal O}_X(P_0)=L,\quad\hbox{and }
{\rm deg}V=d.$$
Such elementary transformations $\varphi$ are classified by surjective maps in
$$Hom({\cal O}_X^{\oplus n}, {\cal O}_P)=\hbox{$\bigoplus\limits_{i=1}^d$}
Hom ({\cal O}_X^{\oplus n},{\cal O}_{p_i})\buildrel\over= ({\Bbb C}^n)^d.
\eqno (3.2)$$
\end{say}

\begin{notationnum} A map $\varphi$ in (3.1)
 can be expressed, under the isomorphisms in (3.2),  as a matrix
$$\left( \begin{array}{c}
\varphi_1 \\ \vdots \\ \varphi_d 
\end{array} \right) =
\left( \begin{array}{ccc}
a_{11}&\ldots&a_{1n} \\ 
\vdots&\ddots&\vdots \\ 
a_{d1}&\ldots&a_{dn} 
\end{array} \right) \eqno (3.3) $$
where $\varphi_i$ represents the vector $(a_{i1},\ldots, a_{in})$ in
$Hom({\cal O}_X^{\oplus n}, {\cal O}_{p_i}) \cong {\Bbb C}^n$. 
\end{notationnum}

\begin{rem} 
\label{norowzero}
$\varphi$ is a surjection if and only if no rows 
$\varphi_i$ in (3.3) are zero vectors. 
\end{rem}

\begin{say}
Tensor  the exact sequence (3.1) by $K_X$ and then take the long 
cohomological exact sequence, we get
$$0\longrightarrow H^0(V^*\otimes K_X)\longrightarrow H^0({\cal O}_X^{\oplus n}
\otimes K_X)\mapright {\varphi^0} H^0({\cal O}_P) $$
$$\longrightarrow H^1(V^*\otimes K_X)
\longrightarrow  H^1({\cal O}_X^{\oplus n}\otimes K_X)\longrightarrow 0. \eqno (3.4) $$
\end{say}

 We need to introduce two generic conditions on elementary transformations (3.1). 
We shall then show  that if an elementary transformation satisfies these two 
conditions, then the vector bundle $V$ is stable and  $h^1(V)=0$.

\begin{condition} We define the following:
\begin{enumerate}
\item{(A)} An elementary transformation  $\varphi$ in 
$Hom({\cal O}_X^{\oplus n},{\cal O}_P)$ is said to satisfy
Condition A if the induced map $\varphi^0$ of $\varphi$ in (3.4)
$$\varphi^0\colon H^0(K_X^{\oplus n})\mapright{} H^0({\cal O}_P)$$
 is an isomorphism.
\item{(B)} An elementary  transformation  $\varphi$ in 
$Hom({\cal O}_X^{\oplus n},{\cal O}_P)$ is said to satisfy
Condition B if  any $n$ many $\varphi_i$'s
in  (3.3) are linearly independent.
\end{enumerate}
\end{condition}

\begin{rem}
\label{conditionsab}
\par
\begin{enumerate}

\item From the exact sequence (3.4), we see that if $\varphi$ satisfies
Condition A, then $h^0(V^*\otimes K_X)=0$, or $h^1(V)=0$ by Serre duality. 
Converse is also true. That is,  if $h^1(V)=0$, then $\varphi^0$ is injective.
Because $h^0({\cal O}_X^{\oplus n} \otimes K_X) = h^0({\cal O}_P) = d$, it has to be 
an isomorphism. Hence Condition A is equivalent to $h^1(V)=0$. By Riemann-Roch, it is also
equivalent to $h^0(V)=n$.
\item Generic $d\times n$ matrix satisfies Condition B.  
 Also this property is clearly invariant under the natural 
action of  $Aut({\cal O}_X^{\oplus n})$ and $Aut ({\cal O}_P)$
on $Hom({\cal O}_X^{\oplus n}, {\cal O}_P)$.
\end{enumerate}
\end{rem}

\begin{lem} Fix any divisor $P$ in $U$.  Generic  elementary transformation
$\varphi$ in  (3.1) satisfies  Condition A.
\end{lem}
\proof  First  notice that $h^0(K_X^{\oplus n})=ng =d =h^0({\cal O}_P)$. 
Let $\omega_1,\ldots,\omega_g$ be a basis of $H^0(K_X)$ as chosen in 
\ref{choosedivisor}.
Take
$$\underbrace{(\omega_1, 0,\ldots,0)}_n,\ldots,(\omega_g,0\ldots,0),\ldots,
(0,\ldots, 0,\omega_1),\ldots, (0,\ldots, 0, \omega_g)$$
 as a basis of 
$H^0(K_X^{\oplus n})$. One checks that
$$\varphi^0(\omega_1,0\ldots, 0)= (a_{11}\omega_1(p_1),a_{21}\omega_1(p_2),
\ldots, a_{d1}\omega_1(p_d)).$$
\par
By the similar calculation for other elements of the basis, we get a natural 
matrix representation of the map $\varphi^0$:
$$\left( \begin{array}{cccc}
a_{11}\omega_1(p_1)&a_{21}\omega_1(p_2)&\cdots&a_{d1}\omega_1(p_d) \\ 
a_{11}\omega_2(p_1)&a_{21}\omega_2(p_2)&\cdots&a_{d1}\omega_2(p_d) \\ 
\vdots&\vdots&\ddots&\vdots \\ 
a_{11}\omega_g(p_1)&a_{21}\omega_g(p_2)&\cdots&a_{d1}\omega_g(p_d) \\ 
a_{12}\omega_1(p_1)&a_{22}\omega_1(p_2)&\cdots&a_{d2}\omega_1(p_d) \\ 
\vdots&\vdots&\ddots&\vdots \\ 
a_{12}\omega_g(p_1)&a_{22}\omega_g(p_2)&\cdots&a_{d2}\omega_g(p_d) \\ 
\vdots&\vdots&\ddots&\vdots \\ 
a_{1n}\omega_1(p_1)&a_{2n}\omega_1(p_2)&\cdots&a_{dn}\omega_1(p_d) \\ 
\vdots&\vdots&\ddots&\vdots \\ 
a_{1n}\omega_g(p_1)&a_{2n}\omega_g(p_2)&\cdots&a_{dn}\omega_g(p_d) 
\end{array} \right) \eqno (3.5)$$
Hence  $\varphi$ satisfies  Condition A iff the determinant
of (3.5) is not zero.

Regarding  the determinant of (3.5) as a polynomial in variables $a_{ij}$ with
$\omega_i(p_j) $'s fixed,  either the polynomial is identically zero or the
zero locus of this polynomial is a divisor in the space $Hom({\cal O}_X^{\oplus
n}, {\cal O}_P)$. 
However we know that if $P$ is in $U$, the determinant of  the matrix (3.5) 
is not  zero at $a_{ij}=b_{ij}$ (see \ref{chooseU}).
Hence the set of the  elementary transformations $\varphi$ satisfying
Condition  A is a non-empty  open dense
subset of the space $Hom({\cal O}_X^{\oplus n}, {\cal O}_P)$.\endproof

\begin{cor} For generic elementary 
transformations $\varphi\in Hom({\cal O}_X^{\oplus n}, {\cal O}_P)$,
$$0\longrightarrow V^*\mapright{}{\cal O}_X^{\oplus n}\mapright{\varphi}
{\cal O}_P\mapright{}0 \qquad \hbox{where $P\in U$},$$                        
$h^0(V)=n$ or equivalently
$h^1(V)=0$.
\end{cor}

\begin{thm} 
\label{stable}
Fix a divisor $P\in U$. 
If  an elementary transformation   $\varphi$  satisfies Condition A and B, then
 $V$ is stable and  satisfies:
\begin{enumerate}
 \item  ${\rm deg} V =ng$ and ${\rm det}V={\cal O}_X(P)=L$;
\item  $h^1(V)=0 $.
\end{enumerate}
\end{thm}
\proof  Suppose $V$ is not stable. Then there exist rank-$f$ vector
bundle $F$, rank-$f'$ vector bundle $F'$,  and an exact sequence:
$$\exact{F}{V}{F'}\eqno(3.6)$$
such that ${\rm deg}F\ge \displaystyle {f\cdot{\rm deg} V\over n} =fg$.

Take the dual of (3.1), we get an exact sequence
$$\exact{{\cal O}_X^{\oplus n}}{V}{{\cal O}_P}.\eqno(3.7)$$

Combining (3.6) and (3.7), we obtain a commutative diagram of exact sequences:
$$\begin{array}{cccccccccccc}
&&&&0 \\ 
&&&&\mapup{} \\ 
&&&&F' \\ 
&&&&\mapup{} \\ 
0&\mapright{}&{\cal O}_X^{\oplus n}&\mapright{}&V&\mapright{}&{\cal O}_P&
\mapright{}&0 \\ 
&&\mapup{}&&\mapup{}&&\mapup{} \\ 
0&\mapright{}&E&\mapright{}&F&\mapright{}&{\cal O}_Q&\mapright{}&0 \\ 
&&\mapup{}&&\mapup{}&&\mapup{} \\ 
&&0&&0&&0
\end{array} \eqno (3.8) $$
where $E$ is a rank-$f$ vector bundle and $Q$ is a subset of $P$.
\par
Diagram (3.8) can be extended into another commutative diagram of 
exact sequences:

$$\begin{array}{ccccccccccc}
&&0&&0&&0 \\ 
&&\mapup{}&&\mapup{}&&\mapup{} \\ 
0&\mapright{}&E'&\mapright{}&F'&\mapright{}&{\cal O}_R&\mapright{}&0 \\ 
&&\mapup{}&&\mapup{}&&\mapup{} \\ 
0&\mapright{}&{\cal O}_X^{\oplus n}&\mapright{}&V&\mapright{}&{\cal O}_P&\mapright{}&0 \\ 
&&\mapup{}&&\mapup{}&&\mapup{} \\ 
0&\mapright{}&E&\mapright{}&F&\mapright{}&{\cal O}_Q&\mapright{}&0 \\ 
&&\mapup{}&&\mapup{}&&\mapup{} \\ 
&&0&&0&&0
\end{array} \eqno (3.9)$$
where $R$ is a subset of $P$.

>From Riemann-Roch, we get $h^0(F)={\rm deg}F+f(1-g)+h^1(F)\ge f$ and equality
holds iff $h^1(F)=0$ and ${\rm deg} F=fg$.
\par
Take long cohomological exact sequence of the diagram (3.9), we get
a commutative diagram of exact sequences:
$$\begin{array}{ccccccccccccccccc}
0&\mapright{}&H^0(E')&\mapright{\alpha_1}&H^0(F')&\mapright{\alpha_2}&
H^0({\cal O}_R) \\ 
&&\mapup{}&&\mapup{\phi_1}&&\mapup{\varphi_1} \\ 
0&\mapright{}&H^0({\cal O}_X^{\oplus n})&\mapright{\beta_1}&
H^0(V)&\mapright{\beta_2}&H^0({\cal O}_P) \\ 
&&\mapup{}&&\mapup{\phi_2}&&\mapup{\varphi_2} \\ 
0&\mapright{}&H^0(E)&\mapright{\gamma_1}&H^0(F)&\mapright{\gamma_2}
&H^0({\cal O}_Q) \\ 
&&\mapup{}&&\mapup{}&&\mapup{} \\ 
&&0&&0&&0. 
\end{array} \eqno (3.10)$$
Since the map $\beta_2$ is a zero map by Condition A, we get $\gamma_2$ is
 also a zero map, hence $h^0(E)=h^0(F)$. 
\par
Assume that $h^0(F)>f$,  then $h^0(E)=h^0(F)>f$. 
Consider the sub-sheaf $E_1$ of $E$ generated by global sections of $E$. Since 
$E$ is a sub-sheaf of ${\cal O}_X^{\oplus n}$, $H^0(E)$ is a subspace of
$ H^0({\cal O}_X^{\oplus n})$. Hence the sub-sheaf $E_1$ is a trivial
sheaf with  rank 
equal to $h^0(E)>f$. This contradicts to the
fact that ${\rm rank}E=f$.

Thus $h^0(F)$ must equal $f$, therefore $h^1(F)=0$ and ${\rm deg} F=fg$.
In this case, we get 
$${\cal O}^{\oplus f}\hookrightarrow E\hookrightarrow {\cal O}^{\oplus n}.$$

 The first map forces ${\rm deg} E\ge 0$ and the second map forces 
${\rm deg}E\le 0$. Hence ${\rm deg}E=0$ and the map ${\cal O}_X^{\oplus f}
\longrightarrow E$ is an isomorphism. From the exact sequence
$\exact{{\cal O}_X^{\oplus f}}{F}{{\cal O}_Q}$, we can see that the number of points
in $Q$ is $\ell(Q)={\rm deg} F=fg$. Hence $\ell(R)=ng -fg =f'g$.
\par
Now we consider the dual diagram of (3.9):
$$\begin{array}{ccccccccccccc}
{}&&0&&0&&0 \\ 
&&\mapdown{}&&\mapdown {}&&\mapdown{} \\ 
0&\mapleft{}&{\cal O}_R&\mapleft{}&E^{'*}&\mapleft{}&F^{'*}&\mapleft{}&0 \\ 
&&\mapdown{}&&\mapdown {}&&\mapdown{} \\ 
0&\mapleft{}&{\cal O}_P&\mapleft{}&{\cal O}_X^{\oplus n}&\mapleft{}&V^*&
\mapleft{}&0 \\ 
&&\mapdown{}&&\mapdown{}&&\mapdown{} \\ 
0&\mapleft{}&{\cal O}_Q&\mapleft{}&{\cal O}^{\oplus f}_X&\mapleft{}&F^*&
\mapleft{}&0 \\ 
&&\mapdown{}&&\mapdown{}&&\mapdown{} \\ 
&&0&&0&&0. 
\end{array}  \eqno (3.11)$$
\par
Since $h^1(V)=0$ (due to Condition  A),  $h^0(V^*\otimes K_X)=0$ by the Serre duality. 
Also $h^0(K_X)\ge 1$ because we have assumed
that $g\ge 2$. Hence
 $h^0(V^*)=0$ and $h^0(F^{'*})=0$ as well. Since
$h^1(F)=0$, by the similar argument, $h^0(F^*)=0$.
Now we take the cohomology of the diagram (3.11), we get
$$\begin{array}{ccccccccc}
0&&0 \\ 
\mapdown{}&&\mapdown{} \\ 
H^0({\cal O}_R)&\mapleft{\rho_2}&H^0(E^{'*})&\mapleft{}&0 \\ 
\mapdown{\pi_2}&&\mapdown{\sigma_2} \\ 
H^0({\cal O}_P)&\mapleft{\varphi_0}&H^0({\cal O}_X^{\oplus n})&\mapleft{}&0 \\ 
\mapdown{\pi_1}&&\mapdown{\sigma_1} \\ 
H^0({\cal O}_Q)&\mapleft{\rho_1}&H^0({\cal O}_X^{\oplus f})&\mapleft{}&0 \\ 
\mapdown{}&&\mapdown{} \\ 
0&&0 
\end{array}$$
where $\varphi_0$ is the induced map of $\varphi$ on the cohomology. (Notice that the map 
$\varphi_0$ here is different from the map $\varphi^0$ we defined in (3.4).) 
\par
The surjectivity of $\sigma_1$ is due to the nature of the map
 ${\cal O}_X^{\oplus n}\longrightarrow {\cal O}_X^{\oplus f}$.
Hence $h^0(E^{'*})=n-f=f'$.
(In fact, we can show that $E^{'*}={\cal O}_X^{\oplus f'}$, but we don't 
need such a strong  statement.)
\par
$\pi_1$ is a natural coordinate projection. By Condition B, 
$\pi_1|_{{\rm Im}\varphi_0
}$ is either a surjection if $\ell(Q)\le n$ or an injection if $\ell (Q)\ge n$.
\par
If $\ell(Q)\le n$, $\pi_1|_{{\rm Im}\varphi_0}$ is a surjection. Hence $\rho_1$ has to 
be a surjection. But $h^0({\cal O}_X^{\oplus f})=f<fg=h^0({\cal O}_Q)$, 
a contradiction. 
\par
If $\ell(Q)\ge n$, $\pi_1|_{{\rm Im}\varphi_0}$ is an injection. 
${\rm Im}(\pi_1\circ \varphi_0
)$ has rank $n$. But ${\rm Im}(\pi_1\circ \varphi_0)={\rm Im}(\rho_1\circ\sigma_1)$ and
${\rm Im}(\rho_1\circ\sigma_1)$ has rank $f<n$, a contradiction.
\par
Hence $V$ has to be stable.
\endproof

%\begin{rem} Grzegorczyk
% also proved (2) in his paper [{\bf G}].
%\end{rem}

\section{Group actions on $Hom({\cal O}_X^{\oplus n}, {\cal O}_P)$}

 In this section, we study the  actions of $Aut ({\cal O}_X^{\oplus n})$ and $Aut ({\cal O}_P)$
on $Hom({\cal O}_X^{\oplus n}, {\cal O}_P)$.

\begin{say} Recall that every surjective map $\varphi$ in $Hom({\cal O}_X^{\oplus n},
{\cal O}_P)$ gives an exact sequence 
$$0\longrightarrow {V^*}\mapright{\phi}
 {{\cal O}_X^{\oplus n}}\mapright{\varphi}
{\cal O}_P\longrightarrow 0.$$
Hence for an elementary transformation $\varphi$, we get a vector bundle
$V^*$, or equivalent $V$. The group $Aut({\cal O}_P)$ and
 $Aut({\cal O}_X^{\oplus n})$ act on $Hom({\cal O}_X^{\oplus n}, {\cal O}_P)$
as follows. Let  $\rho\in Aut ({\cal O}_X^{\oplus n})$ 
and $\sigma\in Aut({\cal O}_P)$. The action of $\rho$ and $\sigma$ on 
$\varphi$ gives a new 
elementary transformation $\varphi'=\sigma\circ \varphi\circ \rho^{-1}$.
It fits into the following exact sequence:
$$0\longrightarrow V^* \buildrel{\rho\circ \phi}\over
\longrightarrow {\cal O}_X^{\oplus n}\mapright
{\varphi'}{\cal O}_P\longrightarrow 0.$$
This means that $\varphi'$ gives arise to the same bundle $V^*$, or equivalently, $V$. So
we have:
\end{say}

\begin{lem}
\label{orbitimpliesisom}
 Suppose two elementary transformations $\varphi,
\varphi'\in Hom({\cal O}_X^{\oplus n}, {\cal O}_P)$ are in the same orbit
under the actions of $Aut ({\cal O}_X^{\oplus n})$ and 
$Aut ({\cal O}_P)$. Then the corresponding vector bundles (which are kernels of them)
 are isomorphic.
\end{lem}

 It is natural to ask whether the converse is also true. In general, it
might well be that two elementary transformations in different group orbits 
give arise to two isomorphic vector bundles. However we have the following
lemma:

\begin{lem} 
\label{isomimpliesorbit}
Let $V$ and $V'$ be obtained from two 
elementary transformations $\varphi$ 
and  $\varphi'	$ in $Hom({\cal O}_X^{\oplus n}, {\cal O}_P)$
and  $Hom({\cal O}_X^{\oplus n}, {\cal O}_Q)$
respectively  both of which satisfy Condition  A ($P$
and $Q$ are not assumed to be the same). If $V$ and $V'$ are isomorphic, 
then $P=Q$ and $\varphi$ and $\varphi'$ are in the same
group orbit under the  actions of $Aut({\cal O}_X^{\oplus n})$ and 
$Aut ({\cal O}_P)$. 
\end{lem}
\proof Let $f\colon V\mapright{} V'$ be the isomorphism. 
 Consider the following exact sequences 
$$0\longrightarrow V^*\longrightarrow {\cal O}_X^{\oplus n}\mapright{\varphi}
{\cal O}_P\mapright{} 0, $$
$$0\longrightarrow V^{'*}\longrightarrow {\cal O}_X^{\oplus n}\mapright{\varphi'}
{\cal O}_Q\mapright{} 0. $$
Take the duals of these two exact sequences, we have 
$$0\mapright{}{\cal O}_X^{\oplus n}\mapright{\alpha}V\mapright{\beta}
{\cal O}_P\mapright{}0,\eqno(4.1)$$
$$0\mapright{}{\cal O}_X^{\oplus n}\mapright{\alpha'}V'\mapright{\beta'}
{\cal O}_Q\mapright{}0.\eqno(4.2)$$
Since $\varphi$ and $\varphi'$ satisfy Condition  A, from Remark \ref{conditionsab},
 we have that
$h^0(V)=h^0(V')= n$. Take  the 
cohomology of the two exact sequences (4.1) and (4.2), 
we get two isomorphisms
$$H^0({\cal O}_X^{\oplus n})\mapright{\alpha^0} H^0(V),\qquad
H^0({\cal O}_X^{\oplus n})\mapright{\alpha^{\prime 0}} H^0(V')$$
where $\alpha^0$ (or $\alpha^{'0}$) is the induced map of $\alpha$
 (or $\alpha^{'}$) on
cohomologies. 
Hence ${\rm Im}\alpha$ (or ${\rm Im}\alpha'$) 
 is the sub-sheaf of $V$ (or $V'$) generated by $H^0(V)$ (or $H^0(V')$).
Since $H^0(V)$ is mapped isomorphically to $H^0(V')$ by the map 
$f^0$ which is the induced map of the given isomorphism $f: V \rightarrow V' $,
 ${\rm Im}\alpha$ is mapped isomorphically to ${\rm Im}\alpha'$ by the 
map $f$ from $V$ to $V'$. Therefore, there is an induced morphism 
$\rho\in Aut({\cal O}_X^{\oplus n})$  
making the following diagram commute:
$$\begin{array}{cccccccccccc}
0 \mapright{}&{\cal O}_X^{\oplus n}&\mapright{\alpha}&V \\ 
&\mapdown{\rho}&&\mapdown{f} \\ 
0\mapright{}&{\cal O}_X^{\oplus n}&\mapright{\alpha'}&V' 
\end{array} $$
which in turn induces a morphism $\sigma\colon
{\cal O}_P\mapright{}{\cal O}_Q$. That is to say, 
there exists a commutative diagram:
$$\begin{array}{cccccccccccc}
0\mapright{}&{\cal O}_X^{\oplus n}&\mapright{\alpha}&V&
\mapright{\beta}&{\cal O}_P&\mapright{}0 \\ 
&\mapdown{\rho}&&\mapdown{f}&&\mapdown{\sigma} \\ 
0\mapright{}&{\cal O}_X^{\oplus n}&\mapright{\alpha'}&
V'&\mapright{\beta'}& {\cal O}_Q&\mapright{}0. 
\end{array} \eqno(4.3)$$

Notice that points in $P$ (or $Q$) 
are the points where the map $\alpha$ (or $\alpha'$) 
fails to be an isomorphism. 
Hence $Q$ must equal $P$ and $\sigma$ must be an automorphism of ${\cal O}_P$.
\par
Take the dual of (4.3), we get 
$$\begin{array}{cccccccccccc}
0\mapright{}&V^*&\mapright{}&{\cal O}_X^{\oplus n}&
\mapright{\varphi}&{\cal O}_P&\mapright{}0 \\ 
&\mapup{f^*}&&\mapup{\rho^*}&&\mapup{\sigma^*} \\ 
0\mapright{}&V^{'*}&\mapright{}&{\cal O}_X^{\oplus n}&
\mapright{\varphi'}&
{\cal O}_P&\mapright{}0 
\end{array}$$
where $\rho^*\in Aut ({\cal O}_X^{\oplus n})$ and $\sigma^*\in 
Aut ({\cal O}_P)$. So $\varphi =\sigma^*\circ \varphi'\circ (\rho^*)^{-1}$.
\par
Hence $\varphi$ and $\varphi'$ are in the same group orbit.
\endproof

\begin{say}
Lemma \ref{orbitimpliesisom} and Lemma \ref{isomimpliesorbit} 
tell us that for generic elementary transformation
$\varphi\in Hom({\cal O}_X^{\oplus n}, {\cal O}_P)$, its group orbit
$Aut({\cal O}_P)\cdot \varphi\cdot Aut({\cal O}_X^{\oplus n})$
classifies stable vector bundles up to bundle isomorphisms. 
This leads us to the quotient space $$
 Aut ({\cal O}_P) \backslash Hom({\cal O}_X^{\oplus n}, {\cal O}_P)/
Aut({\cal O}_X^{\oplus n}).$$
\end{say}

\begin{say}
We have the following identifications:
$$Hom({\cal O}_X^{\oplus n},{\cal O}_P)\buildrel \sim\over=
H^0({\cal O}_P\otimes {\cal O}_X^{\oplus n}) $$
$$\cong H^0({\cal O}_{p_1}^{\oplus n})\oplus \ldots\oplus
H^0({\cal O}_{p_d}^{\oplus n}) \cong ({\Bbb C}^n)^d. \eqno(4.4)$$
Under the isomorphism in (4.4), the group actions  become the natural ones:
\begin{enumerate}
\item  $Aut ({\cal O}_X^{\oplus n})=GL(n)$ 
acts on the  space 
$\bigoplus\limits_{i=1}^d H^0({\cal O}_{p_i}^{\oplus n})=({\Bbb C}^n)^d$ 
diagonally (on each component $H^0({\cal O}_{p_i}^{\oplus n})={\Bbb C}^n$
the action is the standard one);
\item  $Aut ({\cal O}_P^{\oplus n})= ({\Bbb C^*})^d$ acts on  
$Hom({\cal O}_X^{\oplus n}, {\cal O}_P)=({\Bbb C}^n)^d$ by component-wise scalar multiplications.
\end{enumerate}
In a more down-to-earth way,  we may write an element $\varphi\in ({\Bbb C}^n)^d$ as
 $(\varphi_1,
\ldots, \varphi_d)$ where each $\varphi_i\in {\Bbb C}^n$. Then for $g\in GL(n)$,
$$g\cdot \varphi =(g\varphi_1,\ldots, g\varphi_d),$$
and for $c=(c_1,\ldots, c_d)\in  ({\Bbb C^*})^d$,
$$c\cdot \varphi=(c_1\varphi_1,\ldots, c_d\varphi_d).$$
\end{say}

\begin{lem}
\label{invariant}  Condition A (or B) is invariant under
the actions of $Aut({\cal O}_P)$ and $Aut({\cal O}_X^{\oplus n})$. That is, if
$\varphi$ is an elementary transformation in 
$Hom({\cal O}_X^{\oplus n}, {\cal O}_P)$ satisfying  Condition A (or
B), then every element in the orbit $Aut({\cal O}_P)\cdot \varphi\cdot
Aut({\cal O}_X^{\oplus n})$ also satisfies Condition  A (or B). 
\end{lem}

\proof From Remark \ref{conditionsab}, we know that Condition A is equivalent to
 $h^1(V)=0$. From Lemma \ref{orbitimpliesisom}, we know that $\varphi$ and 
 $\sigma\circ\varphi\circ \rho^{-1}$ 
 induce the same vector bundles $V$ for any $\sigma\in
Aut({\cal O}_P)$ and $\rho\in Aut ({\cal O}_X^{\oplus n})$. Thus
if $\varphi$ satisfies Condition  A, so does  $\sigma\circ\varphi\circ \rho^{-1}$.

The invariance of Condition B under the actions of $Aut({\cal O}_P)$ and 
$Aut ({\cal O}_X^{\oplus n})$ is straightforward from the definition of 
Condition B. 
\endproof

\begin{notationnum}

\begin{enumerate}
\item  Define $N$ to be the space of elementary transformations in 
$Hom({\cal O}_X^{\oplus n}, {\cal O}_P)$, i.e. $$N\buildrel\sim\over=
({\Bbb C}^n-\{{\bf 0}\})^d$$
where ${\bf 0}$ is the zero vector in ${\Bbb C}^n$.
\item Define $\overline{N}$ to be the quotient space of $N$ by the 
action of $Aut({\cal O}_P)$, i.e.
$$\overline{N}= ({\Bbb C}^*)^d\backslash ({\Bbb C}^n-\{{\bf 0}\})^d
 \buildrel\sim\over= ({\Bbb P}_{n-1})^d.$$
\item Define $\widehat{N}$ to be the (non-separated) quotient space of
$N$ by the action of groups $Aut({\cal O}_P)$ and $Aut({\cal O}_X^{\oplus n})$,
i.e.
$$\widehat{N} \buildrel\sim\over=({\Bbb C}^*)^d\backslash
({\Bbb C}^n-\{{\bf 0}\})^d/GL(n)
\buildrel\sim\over=({\Bbb P}_{n-1})^d/PGL(n)$$
where $PGL(n)$ acts on $({\Bbb P}_{n-1})^d$ diagonally.
\item Define $N_A$ (or $N_B$) to be the subset of elementary transformations
in $N$ satisfying Condition  A (or B); define $\overline{N}_A$ (or $\overline{N}_B$)
to be the quotient space of $N_A$ (or $N_B$) by the action of the group
$Aut({\cal O}_P)$; and define $\widehat{N}_A$ (or $\widehat{N}_B$) to be the quotient 
space of $N_A$ (or $N_B$) by the action of the groups $Aut({\cal O}_P)$ and $Aut({\cal O}_X^{\oplus n})$.
Clearly $\widehat{N}_A$ and $\widehat{N}_B$ are non-empty
Zariski open subsets of $\widehat{N}$ and 
we shall see that $\widehat{N}_B$ is separated, quasi-projective, and rational (see \ref{git} below).
\end{enumerate}
\end{notationnum}

\begin{say}
\label{git}
Now take an element 
 $\varphi=(\varphi_1,\ldots,\varphi_d)$ 
in $N$. Consider its image 
$\overline\varphi=(\overline\varphi_1,\ldots,\overline\varphi_d)$ in $\overline
N$. $PGL(n)$ acts on $\overline N=({\Bbb P}_{n-1})^d$.
 The geometric invariant theory of this standard $PGL(n)$ action on
 $\overline N=({\Bbb P}_{n-1})^d$ can be found in [{\bf MF}]. If $\varphi$
satisfies Condition  B, then any $n$-many  
$\overline\varphi_{i_1},\ldots,\overline\varphi_{i_n}$ will be linearly
independent. In the context of geometric invariant theory, 
 such $\overline\varphi$'s are necessarily stable with respect to all linearizations.
It follows from [{\bf MF}] that $\widehat{N}_B$ is separated, quasi-projective, and rational.

For the convenience to the reader, we shall give a brief account of these standard results.
Consider the diagonal action of $PGL(n+1)$ on $({\Bbb P}_n)^{m+1}$.
Let $p_i$ be the projection of $({\Bbb P}_n)^{m+1}$ to the $i$-th factor. 
Let $L_i=p_i^*{\cal O}_{{\Bbb P}_n}(1)$. Then ${\cal L}=(L_1\otimes \ldots
\otimes L_{m+1})^{n+1}$ admits a $PGL(n+1)$-linearization.
%(more correctly,
%we should replace $PGL(n+1)$-action by $GL(n+1)$-action and consider
%$GL(n+1)$-linearizations instead because $PGL(2)$ may not admit any linearization.
%The quotients, however, will not be affected). 
 Assume that $m\ge n+1$.

The following is the Definition 3.7 / Proposition 3.4 in [{\bf MF}]. 
\begin{prop}
 The set of stable points
in $({\Bbb P}_n)^{m+1}$ with respect to $\cal L$ is the open subset of 
$({\Bbb P}_n)^{m+1}$ whose geometric points $x=(x^{(0)}, \ldots, x^{(m)})$ are 
those points such that for every proper linear subspaces $L\subset {\Bbb P}_n$,
$${\hbox{ number of points $x^{(i)}$ in $L$}\over m+1}<{ {\rm dim} L+1\over n+1.}
\eqno(4.5)$$
\end{prop}

It is an easy exercise to check the following:

\begin{cor} Let  $x=(x^{(0)}, \ldots, x^{(m)})$ be a
point in $({\Bbb P}_n)^{m+1}$ such that any $(n+1)$-many $x^{(i)}$ are linearly
independent in ${\Bbb P}_n$, then $x$ is a stable point
 with respect to $\cal L$.
\end{cor}

\begin{defn}
\label{defnstable}
 We say an elementary transformation 
$\varphi$ in $Hom({\cal O}_X^{\oplus n}, {\cal O}_P)$ is stable if 
$\overline \varphi\in\overline N$ is  stable with respect to the linearization 
$\cal L$. 
\end{defn}
\end{say}
Apply the geometric invariant theory to our situation, we get:
\begin{cor} An elementary transformation $\varphi\in Hom({\cal O}_X^{
\oplus n}, {\cal O}_P)$  is stable if $\varphi$
satisfies Condition B.
\end{cor}

Now we may rewrite Theorem \ref{stable} as follows:

\begin{prop}
\label{re-stable} Fix a divisor $P$ in $U$. Let $\varphi$ be an elementary
transformation in $Hom({\cal O}_X^{\oplus n}, {\cal O}_P)$. 
If $\varphi$ satisfies Condition  A and B, then the 
vector bundle $V$ obtained from $\varphi$ is a stable vector bundle with 
${\rm deg}V=d=ng$ and ${\rm rank} V=n$ such that
\begin{enumerate}
\item $h^1(V)=0$.
\item ${\varphi}$ is stable (in the sense of Definition \ref{defnstable}). 
\item The set of such stable bundles is isomorphic to the set
$\widehat{{ N}}_A\cap\widehat{{ N}}_B$. Furthermore, $\widehat{{ N}}_A\cap\widehat{{ N}}_B$ is 
quasi-projective, rational, and  of dimension $d(n-1)-n^2+1$.
\end{enumerate}
\end{prop}

\begin{say}
Now we see that we can construct a stable vector bundle such that it satisfies (1),
 (2) and (3) of Theorem \ref{constructingstablebundles}. Next, we need to prove that such vector
bundles are generic in the  moduli space ${\cal M}(n, L)$.

For this, we have to let $P$ move in $U$. Then we get a set $\cal B$ of 
 stable vector bundles. $\cal B$ has dimension
$${\rm dim} \widehat{{ N}}_A\cap\widehat{{ N}}_B+{\rm dim} U
=(n-1)ng-n^2+1+(n-1)g=(n^2-1)(g-1)$$
which is the same as the dimension of the moduli space ${\cal M}(n, L)$.
It is natural to ask if this space $\cal B$ is actually a non-empty
Zariski open subset of the moduli space.
Indeed, we shall show that this is the case in the next section.
\end{say}

\section{A birational model for ${\cal M}(n, L)$}

In this section, using relative extension sheaf, we will
construct a family of stable bundles which provides an injection to
the moduli space ${\cal M}(n, L)$ and gives a non-empty Zariski
open subset of ${\cal M}(n, L)$.

In the previous sections, the method we used to construct stable vector bundles is
the elementary transformations. Here we shall use  extensions instead.
The two approaches are related: one is the dual of the other. Below we will elaborate
on this.

\begin{say}
 Given  an elementary transformation $\varphi\in Hom({\cal O}_X^{\oplus n},
{\cal O}_P)$, i.e., 
$$0\longrightarrow {V^*}\longrightarrow
{{\cal O}_X^{\oplus n}}\mapright{\varphi}{{\cal O}_P}\longrightarrow 0\eqno(5.1)$$
Take the dual of the exact sequence (5.1), we get
$$\exact{{\cal O}_X^{\oplus n}}{V}{{\cal O}_P}.\eqno(5.2)$$
Hence $V$ is an extension of ${\cal O}_P$ by 
${\cal O}_X^{\oplus n}$.
 Such extensions are classified by the extension group $Ext^1({\cal O}_P,
{\cal O}_X^{\oplus n})$.
There is an isomorphism
$$Ext^1({\cal O}_P, {\cal O}_X^{\oplus n})\buildrel\sim\over =
H^0({\cal E}xt^1({\cal O}_P, {\cal O}_X^{\oplus n})) 
\buildrel\sim\over=H^0({\cal O}_P\otimes {\cal O}_X^{\oplus n})$$ 
$$\buildrel\sim\over=H^0( {\cal O}_X^{\oplus n}|_{p_1})\oplus \ldots 
\oplus H^0({\cal O}_X^{\oplus n}|_{p_d})\buildrel\sim\over = ({\Bbb C}^n)^d. 
$$
\end{say}

The extension (5.2) in general does not give arise to a vector bundle in the 
middle. However, writing an extension class $e$ as $(e_1, \ldots, e_d) \in ({\Bbb C}^n)^d$,
 we have the following:

\begin{lem} The extension $e$ gives  a vector 
bundle $V$ in (5.2) if and only iff $e_i\ne 0$ for all $i$.
\end{lem}

\proof  See Lemma 16 of [{\bf B1}].\endproof

\begin{rem} The above lemma is equivalent to Remark \ref{norowzero}.
Actually there exists a dictionary between elementary transformations and their
corresponding extensions. For example, under this dictionary,
the groups $Aut({\cal O}_P)$ and $Aut ({\cal O}_X^{\oplus n})$ 
act on $Ext^1({\cal O}_P, {\cal O}_X^{\oplus n})$ in the same way as
they act on $Hom({\cal O}_X^{\oplus n}, {\cal O}_P)$ after we identify
$Ext^1({\cal O}_P, {\cal O}_X^{\oplus n})$ with
$Hom({\cal O}_X^{\oplus n}, {\cal O}_P)$.
We leave the existence of the dictionary to the reader. 
\end{rem}

\begin{notationnum}
\begin{enumerate}
\item  Define $Z'\subset X\times U$ to be  the
universal 
divisor
$$Z'=\{(x,D)\in X\times U|x\in D\}.$$
\item  Define $\pi_i$ to be the natural projection from
$X\times  U$ to the $i$-th factor.
\item Define ${\cal E}$ to be the relative extension sheaf
$${\cal E}={\cal E}xt^1_{\pi_2}({\cal O}_{Z'}, 
\pi^*_1{\cal O}_X^{\oplus n}).$$
The fiber of $\cal E$ over a point $P\in U$ is isomorphic to $Ext^1({\cal O}_P,
{\cal O}_X^{\oplus n})$.
\item  The general linear group $GL(n)$ acts on ${\cal E}$ via acting on ${\cal O}_X^{\oplus n}$.
Meanwhile, the group scheme ${\cal T} = (\pi'_2)_* {\cal O}_{Z'}^*$ also acts on
${\cal E}$ where $\pi'_2$ is the restriction of $\pi_2$ to $Z'$. 
This group scheme is a twisted torus. The fiber of ${\cal T}$ over a point $P \in U$ is
$Aut({\cal O}_P^{\oplus n} ) \cong ({\Bbb C}^*)^d $. The fiber-wise
action of  ${\cal T}$ on ${\cal E}$ is just the 
action of $Aut({\cal O}_P)$ on $Ext^1({\cal O}_P, {\cal O}_X^{\oplus n})$.
The actions of the group scheme ${\cal T}$ and
the (global) group  $GL(n)$ commute. Hence we get an action of ${\cal G} = {\cal T} \times GL(n)$
where the diagonal multiplicative group ${\Bbb C}^*$ 
appears in both ${\cal T}$ and $GL(n)$. 
\end{enumerate}
\end{notationnum}

\begin{lem}
\label{fiber-wise}  The natural morphism:
$${\cal E}_{[P]}\longrightarrow  
Ext^1({\cal O}_P,{\cal O}_X^{\oplus n})\eqno(5.3)$$
is an isomorphism and
 ${\cal E}$ is a locally free sheaf over $U$.
\end{lem}

\proof Note that dim$Ext^1({\cal O}_P,{\cal O}_X^{\oplus n})=nd$ is a constant.
By Satz 3 of [{\bf BPS}],  $\cal E$ is a locally free
sheaf and the natural morphism (5.3) is an isomorphism. \endproof

\begin{prop}
 Let ${\Bbb V}={\Bbb V}({\cal E}^*)$ be the vector bundle
associated to the locally free sheaf $\cal E$ following Grothendieck's
 notation. Then over ${\Bbb V}$, there is a universal extension
$$\exact{q_1^*{\cal O}_X^{\oplus n}\otimes q_2^*{\cal A}}
       {{\cal V}}{\gamma^*{\cal O}_{Z'}}\eqno(5.4)$$
where $q_i$ is the projection from $X\times {\Bbb V}$ to its $i$-th factor
and $\gamma$ is the projection from $X\times {\Bbb V}$ to $X\times U$ and
${\cal A}$ is some line bundle on ${\Bbb V}$. Moreover, $\cal V$ is flat
over $\Bbb V$. 
\end{prop}

\proof Since $Ext^0({\cal O}_P,{\cal O}_X^{\oplus n})=0$ 
for all $P\in U$, ${\cal E}xt^0_{\pi_2}({\cal O}_{Z'}, 
\pi^*_1{\cal O}_X^{\oplus n})=0$. By Lemma \ref{fiber-wise}, 
${\cal E}xt^1_{\pi_2}({\cal O}_{Z'}, 
\pi^*_1{\cal O}_X^{\oplus n})$ commutes with base change
(according to Lange's terminology [{\bf L}]). Hence
by Corollary 3.4 of [{\bf L}], we get the universal extension (5.4). 

Since ${\cal O}_{Z'}$ is flat over $U$, both $\gamma^*{\cal O}_{Z'}$ and
$q_1^*{\cal O}_X^{\oplus n}\otimes q^*_2{\cal A}$ are flat over ${\Bbb V}$.
Hence ${\cal V}$ is flat over $\Bbb V$.
\endproof

\begin{rem}
If we use $e$ to represent a point in ${\Bbb V}$ corresponding to an extension
(5.2), then the restriction  of (5.4) to $X\times e$ is just the extension
(5.2) and ${\cal V}|_{X\times e}\buildrel\sim\over= V$.
\end{rem}

\begin{thm} 
 \label{birationalmodel}
The moduli space ${\cal M}(n, L)$ is birational to
the quotient space of a non-empty
Zariski open subset  of ${\Bbb V}$ by the action
of the group ${\cal G}={\cal T} \times GL(n)$.
\end{thm}
\proof 
Let ${\Bbb V}^0$ be the subset of ${\Bbb V}$ consisting of extensions
(5.2) whose corresponding elementary transformations (via the dictionary)
satisfy Conditions A and B. By Lemma \ref{invariant},
${\Bbb V}^0$ is invariant under the action  of ${\cal G}$. 
By Proposition \ref{re-stable} and semi-continuity using the flatness of
$\cal V$, ${\Bbb V}^0$ is a non-empty Zariski open subset of 
${\Bbb V}$. Restrict ${\cal V}$ to $X\times {\Bbb V}^0$, we get a
family of stable bundles over $X\times {\Bbb V}^0$. Since the moduli 
space ${\cal M}(n, L)$ is coarse, the family induces a morphism:
$$\Phi\colon {\Bbb V}^0\longrightarrow {\cal M}(n,L).\eqno(5.5)$$
Now by Lemma \ref{orbitimpliesisom} and Lemma \ref{isomimpliesorbit}, 
 we see that points of ${\Bbb V}^0$ are in the same ${\cal G}$
orbit if and only if
they correspond to isomorphic stable bundles. 
Thus by passing to the quotient, we get a natural  induced map
$$\overline\Phi \colon {\Bbb V}^0/({\cal T} \times GL(n)) \longrightarrow  {\cal M}(n, L)$$
which  is an injective  morphism. By calculating dimensions, we get
$${\rm dim}({\Bbb V}^0/{\cal G})={\rm dim}{\Bbb V}^0-{\rm dim} ({\cal T} \times GL(n)) + 1 
={\rm dim}{\cal M}(n, L).$$
(Notice that a ``$+1$'' modification appears in the first equality because  
the multiplicative group ${\Bbb C}^*$ appears in both ${\cal T}$ and $GL(n)$.)
Hence the morphism $\overline\Phi$ is birational.\endproof

\begin{rem}
It can be seen that the non-empty Zariski open subset ${\Bbb V}_0/{\cal G}$
of ${\cal M}(n, L)$ is a fibration over  the rational variety $U$
whose typical fiber  is  also  rational since it is birational to the geometric
quotients of $({\Bbb P}_{n-1})^d$ by $PGL(n)$ 
(which are  known to be rational) (see also [{\bf T1, T2}]).
Unfortunately, this fibration seems not to be locally 
trivial in Zariski topology. The twisted torus  ${\cal T}$ is responsible
for the problem. 
\end{rem}

\begin{say} {\sl proof of Theorem \ref{constructingstablebundles}:}
Theorem \ref{constructingstablebundles} clearly follows from the combination of 
Proposition \ref{re-stable} and Theorem \ref{birationalmodel}.
\end{say}

\section{Some remarks on rationality}

We shall provide a thorough and {\sl systematic} account of the inductive method in
[{\bf N1}] and prove
several theorems that are stronger than Theorem \ref{newsteadrationalitythm}.
The essential ideas are, however,  contained in [{\bf N1}].

Throughout,  $(n; d)$ stands for a pair of
integers; while $\gcd (n, d)$ stands for the greatest common factor of the 
two integers $n$ and $d$.

\begin{say} 
Let $V$ be a stable vector bundle in ${\cal M}(n, L)$ and $(n; d)$ satisfies
$$ n(g-1)<d<ng. \eqno(6.1)$$
By Riemann-Roch, $\chi(V)={\rm deg}V+n(1-g)=d-n(g-1).$ Set $r=\chi(V)$. Then $0<r<n$ by (6.1).

In ([{\bf N1}]) and  ([{\bf G}]), it was shown that for  generic bundles $V$ in ${\cal M}(n, L)$, 
$V$  satisfies $h^1(V)=0$ and there exists an exact sequence
$$\exact{{\cal O}_X^{\oplus r}}{V}{V'} \eqno(6.2)$$
such that $V'$ is  stable. That is $V' \in {\cal M}(n', L)$ where $n'=n-r$.
One can always find a non-negative  integer $k$ such that 
$$ n'(g-1)< d'=d-kn' \le n'g. \eqno (6.3)$$
Then ${\cal M}(n', L) \cong {\cal M}(n', L')$ where 
$L'=L \otimes (M^*)^{\otimes n'}$ for some line bundle $M$ of degree $k$ and ${\rm deg}L' = d'$.

This suggests that there exist a rational map from ${\cal M}(n, L)$ to
${\cal M}(n', L')$.
\end{say}

\begin{defn}  Let $(n;d)$  be a pair of positive
integers satisfying (6.1). Let $n'=ng-d$, $d'=d-kn'$ for some non-negative
integer $k$ so that $n'(g-1)<d'\le n'g$. We call 
$(n';d')$ the reduction of  $(n;d)$, we denote this process of reduction by
$$(n;d)\longrightarrow (n';d').$$
\end{defn}

\begin{say} Sometimes, it is necessary to apply  the same 
 procedure to the
dual bundle $V^*$. In other words, instead of constructing $V$
by an extension (6.2), we may do it for $V^*$.
 Precisely, let again $(n;d)$ satisfy (6.1).
Choose a line bundle $M$ with $\hbox{deg}M=2g-1$. Then
$$n(g-1)<\hbox{deg}(V^*\otimes M)=-d+n(2g-1)<ng.$$
Hence $V^*\otimes M$ is a stable vector bundle in the moduli space 
${\cal M}(n, L^*\otimes M^{\otimes n})$. The pair $(n; {\rm deg}(L^*\otimes
M^{\otimes n}))$ satisfies (6.1). 
We then apply the reduction to  $(n; n(2g-1)-d)$, and obtain $(n';d')$ where
$$n'=d-n(g-1), \; d'=n(2g-1)-d-kn',\; \hbox{and} \;n'(g-1) < d' \le n'g $$
for some non-negative integer $k$.
\end{say}

\begin{defn}
 Let $(n;d)$ be a pair of integers 
satisfying the hypothesis (6.1). Let $n'=d-n(g-1)$ and $d'=n(2g-1)-d-kn'$ for
some non-negative integer $k$ such that $n'(g-1)< d' \le n'g$. We call 
 $(n'; d')$ the dual reduction of  $(n; d)$. We denote this process also by
$$(n;d)\longrightarrow (n';d'). $$
\end{defn}

\begin{prop}
Let $(n;d)$ be a pair satisfying (6.1). Let $(n'; d')$ be the reduction (or
dual reduction) of $(n;d)$. Let ${\cal M}(n, L)$ and  ${\cal M}(n', L')$ 
be the moduli spaces such that ${\rm deg} L=d$, ${\rm deg} L'=d'$ and 
$L=L'\otimes M^{\otimes n'}$ for some line bundle $M$ of non-negative degree $k$ 
 with $h^0(M)\ne 0$. Then there exists a non-empty Zariski open
subset  ${\cal M}^0(n, L)\subset  {\cal M}(n, L)$  and a morphism 
$\Phi:  {\cal M}^0(n, L)\mapright{} {\cal M}(n', L')$ such that the image
of $\Phi$ is a  non-empty Zariski open subset of ${\cal M}(n', L')$. 
\end{prop}

\proof We only prove the proposition when $(n;d) \rightarrow (n'; d')$ is a reduction.
The same arguments work for a dual reduction equally well.

Take  ${\cal M}^0(n, L)$ to be the non-empty Zariski open subset of 
${\cal M}(n, L)$ as in Theorem \ref{gthm}.
Recall that a stable vector bundle $V$ in  ${\cal M}^0(n, L)$ has $h^1(V)=0$ 
and sits in the exact sequence (6.2) where $r=n-n'$ and $V'\in
 {\cal M}(n', L)$. 

Let $F'=V'\otimes(M^*)^{\otimes n'}$. Then $F'$ is a stable bundle in
$ {\cal M}(n', L')$. Define the map $$\Phi\colon  {\cal M}^0(n, L)
\mapright{} {\cal M}(n', L')$$ by 
setting
$$\Phi(V)=F'$$

First of all, the map is well-defined. This amounts to saying that if $V_1$
and $V_2$ are two isomorphic stable bundles in ${\cal M}^0(n, L)$, i.e.
there exists two exact sequences 
$$0\mapright{}{\cal O}_X^{\oplus r}\mapright{\alpha_1}V_1\mapright{\beta_1}
V_1'\mapright{}0,$$
$$0\mapright{}{\cal O}_X^{\oplus r}\mapright{\alpha_2}V_2\mapright{\beta_2}
V_2'\mapright{}0$$
where $h^0(V_1)=h^0(V_2)= r$ and $V'_1$ and $V_2'$ are stable, then
$V_1'$ and $V_2'$ are isomorphic. In fact, since $h^0(V_1)=r$, ${\cal O}_X^{
\oplus r}$ is the sub-sheaf of $V_1$ generated by global sections of $H^0(V_1)$. 
Same is true for $V_2$. Hence if $f$ is an isomorphism from $V_1$ to $V_2$,
$f$, restricted to ${\cal O}_X^{\oplus r}$, maps ${\cal O}_X^{\oplus r}
\subset V_1$ to ${\cal O}_X^{\oplus r}\subset V_2$ isomorphically, hence
$V'_1$ is isomorphic to $V'_2$.

Next we  need to show that $\Phi$ is a 
morphism. Consider
 the construction of  ${\cal M}(n, L)$ via geometric invariant 
theory, we let $\cal Q$ be the Quot scheme,  $\cal U$ be the
universal quotient sheaf
 over $X\times {\cal Q}$ and   ${\cal G}$ be the group such 
that $ {\cal M}(n, L)$ is the GIT quotient of $\cal Q$ by $\cal G$. Let
$\Pi$ be the quotient map from ${\cal Q}$ to  ${\cal M}(n, L)$. Let
${\cal Q}_0\subset
 \cal Q$ be the inverse image of $ {\cal M}^0(n, L)$ under the 
map $\Pi$. Since bundles in  ${\cal M}^0(n, L)$ are stable bundles, the 
pre-image of $V\in  {\cal M}^0(n, L)$ under the map $\Pi$ consists of 
bundles isomorphic to $V$.

Clearly we can define a map $$\widetilde \Phi\colon
{\cal Q}_0\mapright{} {\cal M}(n', L')$$ 
in the same way as we defined $\Phi$. 

Because  ${\cal M}(n', L')$  is a coarse moduli space and since $\cal U$ is 
a universal quotient sheaf
 over $X\times {\cal Q}_0$, we conclude that $\widetilde\Phi$
is a morphism.

Since ${\cal M}^0(n, L)$ is the geometric quotient of ${\cal Q}_0$ by ${\cal G}$
and any orbit (e.g., the fiber $\Pi^{-1}(V)$)  is mapped to a single point $F'$ in 
 ${\cal M}(n', L')$ under the map $\widetilde \Phi$, by the universality of
GIT quotients, the morphism $\widetilde \Phi: {\cal Q}_0\mapright{} {\cal M}(n', L')$
induces
 a morphism on the quotient ${\cal M}^0(n, L) \mapright{} {\cal M}(n', L')$.
Evidently, the induced morphism is just the map $\Phi$ defined earlier.

Finally, we need to show that the image of $\Phi$ contains
 a non-empty Zariski open subset of 
${\cal M}(n', L')$.

Recall that the combination of Lemma 5 and Lemma 6 in [{\bf N1}] says that
if $V'$ is stable, $h^1(V')=0$, ${\rm deg} V'=d$, $n(g-1)<d<ng$ and 
$r=d-n(g-1)$, then 
there exists a stable vector bundle $V$ sitting in the exact 
sequence (6.2) with $h^1(V)=0$.

{}From a result in  [{\bf G}], we know that there exists a non-empty Zariski
open subset  $ {\cal M}_0(n', L')$ of ${\cal M}(n', L')$  consisting
 of $F'$
with $h^1(F')=0$. Hence $h^0(F^{'*}\otimes K_X)=0$ by Serre duality. Since
$h^0(M)\ne 0$, we must have $h^0(F^{'*}\otimes (M^*)^{\otimes n'} \otimes K_X)=0$. Recall
that $V'=F'\otimes M^{\otimes n'}$.  Hence
$h^1(V')=h^1(F'\otimes M^{\otimes n'})=0$ by Serre duality. Now
Lemma 5 and Lemma 6 in [{\bf N1}] imply that
there exists a stable vector bundle $V$ in  ${\cal M}^0(n, L)$ such that
$\Phi(V)=F'$. So  ${\rm Im}\Phi$ contains  ${\cal M}_0(n', L')$. Therefore
it contains a  non-empty  Zariski open subset of   ${\cal M}(n', L')$.

The case when $n'=1$, $d'=g$ needs some special remarks.
According to our convention, $V'=L$ in this case. Since $ n\ge 2$, we have
$${\rm deg} L =d\ge n(g-1)+1\ge 2g-2+1=2g-1.$$
Therefore $h^1(L)=h^0(L^*\otimes K_X)=0$ because ${\rm deg}
(L^*\otimes K_X)\le -2g+1+2g-2= -1$. Hence Lemma 5 and Lemma 6 in [{\bf N1}]
also apply.
\endproof

\begin{say} 
\label{genericfibration}
Now we need to analyze the structure of the rational map
$$\Phi\colon {\cal M}(n, L) ---> {\cal M}(n', L').$$ 

Extensions (6.2) are classified by the extension group:
$$Ext^1(V', {\cal O}_X^{\oplus r})=H^1(V^{'*})^{\oplus r}.$$
The group $Aut ({\cal O}_X^{\oplus r})=GL(r)$ acts on the extension group.
Lemma 1 of [{\bf N1}] showed that the $GL(r)$-orbits
correspond to the equivalent classes of stable vector bundles.

Let $F'\in  {\cal M}_0(n', L')$. Then $\Phi^{-1}(F')$ is a non-empty Zariski open
subset of  
$Ext^1(V', {\cal O}_X^{\oplus r})/GL(r)$ which is known to be
rational (see Lemma 2, [{\bf N1}]).
 The rational map $\Phi$ can be regarded as a fibration (over
a non-empty Zariski open subset of ${\cal M}(n', L')$).
\end{say}

\begin{say}
\label{localtriviality} Next, we need to know when the fibration 
$$\Phi\colon {\cal M}(n, L) ---> {\cal M}(n', L')$$
 is locally trivial in Zariski topology over the range of $\Phi$.
 This admits the affirmative answer if 
the moduli space ${\cal M}(n', L')$ is fine, i.e. 
 if ${\cal M}(n', L')$ admits a universal bundle
(if and only if $\gcd(n', d')=1$).  This is because
the existence of the universal bundle will allow us to get
the locally triviality by using the relative extension sheaf (cf. Lemma 3, [{\bf N1}]).
 If  $\gcd(n', d') \ne 1$, 
we still have the map $\Phi$, 
but it seems very likely that the fibration is no longer locally 
trivial in Zariski topology.
\end{say}

%\begin{rem}
%\label{newsteadmistake}
%Among six lemmas in [{\bf N1}], Lemma 1, 2, 5, 6 do not deal with family of
%bundles, hence they do not need the assumption that $(n, d)=1$. 
%Lemma 3 and 4, however, require
%the existence of  universal families.
%\end{rem}

The following is implicit in [{\bf  N2}], which was used to obtain Theorem 
\ref{newsteadrationalitythm}.

\begin{prop} 
\label{keyinduction}
Let $(n;d) \rightarrow (n';d')$ be a reduction or a dual reduction. 
\begin{enumerate}
\item If $(n';d')$ is coprime, then the fibration ${\cal M}(n, L) ---> {\cal M}(n', L')$ is 
      locally trivial in Zariski topology. Furthermore, ${\cal M}(n, L)$
is birational to ${\Bbb C}^m\times {\cal M}(n,' L')$ where $m=(n^2-n^{'2})(g-1)$.
\item If $(n';d')$ is coprime and ${\cal M}(n', L')$ is rational, then ${\cal M}(n, L)$ is rational.
\end{enumerate}
\end{prop}

\begin{say}
\label{reductionsequence}
An interesting question is whether the above induction method works only
for the pairs stated in Theorem \ref{newsteadrationalitythm}. A closer examination
of reductions and dual reductions tells us that the method can actually be extended to work for a much 
larger  class of moduli spaces, which we now describe.

Given a pair $(n;d)$ satisfying (6.1), if $n>1$,
we can apply a 
reduction or a dual reduction to $(n;d)$ to get another pair $(n';d')$. Three possibilities can
occur:
\begin{enumerate}
\item  $n'=1$. 
\item $n' \ge 2$ and  $d'=n'g$.
\item  $n'\ge 2$ and $n'(g-1)<d-kn'< n'g$ for some integer $k\ge 0$.
\end{enumerate}
In the case that (3) occurs,  we can continue  reductions and dual reductions. 
Keep doing this, the process will eventually
terminate and we get a sequence of reductions and dual reductions:
$$(n;d)\longrightarrow (n_1; d_1)\longrightarrow \ldots \longrightarrow (n_t; d_t)\eqno(6.4)$$
such that 
\begin{enumerate}
\item either $n_t=1, d_t=g$, 
\item or $d_t=n_tg$, $n_t\ge 2$.
\end{enumerate}
\end{say}

\begin{defn}
\label{newsteadpair}  Let  $(n; d)$ be a pair of positive
integers satisfying (6.1). The pair is called a  nice pair if either $n=1$ or 
after successive reductions and dual reductions (6.4), we get $n_t=1$.
\end{defn}

\begin{rem} 

 It can happen that after a sequence of reductions and dual reductions,
$(n;d)$ is reduced to $(n_t;d_t)$ with $n_t \ge 2$ and $d_t=n_tg$ while after another 
sequence of reductions and dual reductions, 
$(n;d)$ is reduced to  $(1; g)$. According to our
definition, this kind of pair is nice.
Take the pair $(7; 8)$ when $g=2$ as an example. After two reductions we get
$$(7; 8) \to (6;8) \to (4;8).$$
However, after a dual reduction, we get
$$(7; 8) \to (1;2).$$
The latter will imply that $\cM (7, L)$ with $\deg L = 8$ is rational (cf. Theorem \ref{newsteadpairthm} below),
 while the former won't.
\end{rem}

Below is a simple but useful observation.
\begin{lem}
\label{newsteadpairlem}
Let $(n; d)$ be a pair satisfying (6.1). Let $(n'; d')$ be the reduction
or dual reduction of $(n; d)$. If $\gcd (n', d')=1$, then $\gcd (n, d)=1$.
\end{lem}
\proof
Assume $(n;d)\longrightarrow (n';d')$ is a reduction:
$$(n'; d') = (ng - d; d- k(ng -d))\eqno(6.5)$$
 for some non-negative integer $k$. 	
\par
If an integer $m$ divides $n$ and $d$, from (6.5), it is clear that $m$
divides  $n'=ng-d$ and $d'=d-k(ng-d)$.  Hence the lemma is proved in this case.
\par
Assume  $(n;d)\longrightarrow (n';d')$ is a dual reduction.
$$(n'; d')=(ng-n(2g-1)+d; n(2g-1)-d-k(ng-n(2g-1)+d))\eqno(6.6)$$
 for some non-negative integer $k$. 	
\par
If an integer $m$ divides $n$ and $d$, from (6.6), it is clear that $m$
divides  $n'$ and $d'$ as well.  Hence the lemma is also proved in this case.
\endproof

\begin{rem}
\label{goodpairrem}
\begin{enumerate}
\item On any forwarded (dual) reduction path, $\gcd$ is a non-decreasing function.
\item All nice pairs are coprime.
\item The converse of Lemma \ref{newsteadpairlem}
      is not true which is exactly and the {\it only} place 
       where the proposition of \S 3 in [{\bf N1}] fails.
\end{enumerate}
\end{rem}

\begin{prop}
\label{rationalityinductioncor}
Let $(n;d)\mapright{}(n_1;d_1)\mapright{}\cdots\mapright{}(n_s; d_s)$ be
a sequence of reductions and dual reductions. Suppose ${\rm gcd}(n_s, d_s)=1$, then
${\cal M}(n, L)$ is birational to ${\Bbb C}^x \times {\cal M}(n_s, L_s)$
where ${\deg}L_s = d_s$ and $x=(n^2-n^2_s)(g-1)$.
\end{prop}

\proof By assumption, we have a sequence of 
dominant  rational maps: 
$${\cal M}(n, L) ---> {\cal M}(n_1, L_1)
---> \ldots ---> {\cal M}(n_s, L_s). $$
Since $(n_s, d_s)=1$, all the above moduli spaces are fine moduli spaces by
Lemma \ref{newsteadpairlem}. Hence
each of the above rational map is generically a locally trivial fibration
in Zariski topology whose typical fiber is rational (cf. \ref{localtriviality}
and Proposition \ref{keyinduction}).
The proposition then follows immediately.
\endproof

When $(n_s; d_s)=(1; g)$, we obtain:

\begin{thm}
\label{newsteadpairthm}
If $(n; d)$ is a nice pair, then the moduli space ${\cal M}(n, L)$
with $d={\rm deg} L$ is rational.
\end{thm}

One checks that all pairs in Theorem \ref{newsteadrationalitythm} are nice pairs.
The following gives an example of a nice pair not contained in Theorem \ref{newsteadrationalitythm}.

\begin{exmp}
 Choose $g=6$, $n=15$, $d=77$. One checks that $(n;d)=(15;77)$ satisfies the inequality (3.1) but
satisfies none of the conditions (1), (2), (3) in Theorem 
\ref{newsteadrationalitythm}. 
Apply reductions twice, we obtain
%$$n_1=ng-d=15\times 6-77=13,\qquad d_1=77,$$
%$$n_1(g-1)=13\times 5=65<d_1=77<78=13\times 6=n_1g.$$
%Apply one more reduction to $(n_1;d_1)=(13:77)$, we get
%$$n_2=n_1g-d_1=13\times 6-77=78-77=1.$$
$$(15;77) \rightarrow (13:77) \rightarrow (1;6).$$
Hence the moduli space ${\cal M}(15, L)$ with ${\rm deg}L=77$ is rational.
\end{exmp}

\begin{say}
In view of Theorem \ref{newsteadpairthm},
it is desirable to find ways of constructing all nice pairs. 
The following two lemmas help to characterize arithmetically the iterated
construction of all such pairs starting from $(1;g)$. 
\end{say}

\begin{lem}
\label{iterationlem1} Fix a genus $g \ge 2$.
\begin{enumerate}
\item  $(n;d) \rightarrow (1; g)$  is a (one-step) reduction if and only if 
       $d=ng-1$. Consequently, $\gcd(d, g)=1$.
\item $(n;d) \rightarrow (1; g)$  is a (one-step) dual reduction if and only if 
      $d=ng-n+1$. Consequently, $\gcd(d+n, g)=1$.
\end{enumerate}
\end{lem}

\proof (1) If $d=ng-1$ for $n\ge 2$, it is easy to see that after one
reduction we get $(1;g)$.

Now suppose $(n;d)\rightarrow (1;g)$ is a reduction. From (6.5), we get
$n=\displaystyle{g+(k+1)\over g}$ and $d=g+k$ for some integer $k\ge 0$.
Hence $k+1$ must be $mg$ for some integer $m\ge 1$, 
$n=m+1\ge 2$ and $d=g+mg-1=ng-1$. It  is clear that $\gcd(d, g)=1$.

Similar arguments prove (2).
\endproof

\begin{rem} Lemma \ref{iterationlem1} is 
equivalent to Theorem \ref{newsteadrationalitythm} (1). 
\end{rem}

Next, we have

\begin{lem}
\label{iterationlem2} 
Fix a genus $g \ge 2$. Let $(n'; d')$ be a nice pair with either $\gcd(d', g)=1$ or 
$\gcd(d'+n', g)=1$. Then,  
\begin{enumerate}
\item   there exists   a pair $(n;d)$ having $(n'; d')$ as a reduction 
         if and only if $\gcd(n', g)=1$. In particular, $\gcd(d, g)=1$.
\item    there exists a pair $(n; d)$ having $(n'; d')$ as a dual reduction 
       if and only if $\gcd(n', g)=1$.  In particular $\gcd(d+n, g)=1$.
\end{enumerate}
\end{lem}

\proof
Consider the equations in (6.5):
$$n'=ng-d,\qquad d'=d-kn'.$$
Add these equations together, we get an equation
$$ng=d'+(k+1)n'.$$

Since $\gcd(d', g)=1$ (the case when $\gcd(d'+n', g)=1$ can be proved similarly),
 this equation has integral solutions for $n$ and $k$
iff $\gcd(n', g)=1$. And if such solutions exist, $d=ng-n'$, so $\gcd(d, g)=1$.

The similar arguments prove (2).
\endproof

\begin{cor}
\label{iterationcor}
If $(n; d)$ is a nice pair, then either $\gcd(d, g)=1$ or $\gcd (d+n, g)=1$.
\end{cor}

\begin{rem} In summary, we have the following.
Start with any $(n;d)$. After  successive reductions and dual
reductions, we get
$$(n;d) \rightarrow  (n_1; d_1) \rightarrow \ldots \rightarrow (n_t;d_t).$$
Correspondingly, we have a sequence of dominant  rational maps: 
$${\cal M}(n, L) ---> {\cal M}(n_1, L_1)
---> \ldots ---> {\cal M}(n_t, L_t) \eqno(6.7)$$
where ${\rm deg} L_i=d_i$ fro all $i$.
The terminator is either $(1; g)$ or has that $n_t | d_t$ and $n_t\ge 2$.

The first case is covered in Theorem \ref{newsteadpairthm}.

For the latter case, the sequence (6.7) ends at the moduli space studied in the previous
sections. In the sequence (6.7), 
it is possible that some fibrations are locally
trivial in Zariski topology but the rest are not.  The 
local triviality is known if the base moduli space of a fibration is fine.
Otherwise, we believe that it is not. It seems that
the rationality problem boils down to how well we know about the moduli spaces
${\cal M}(n_t, L_t)$ when  $n_t | d_t$. 
\end{rem}

\begin{say}
The above procedure of proving  rationality (for coprime pairs only)
can be slightly generalized to include an even larger class of pairs 
(essentially due to an observation by Ballico in [{\bf Ba}]). 
The observation rests on the following easy lemma:
\end{say}

\begin{lem} {\rm (Transferring Lemma).}
\label{transferringlemma}
Suppose that we have the following projections:
$$\begin{array}{clc}
X && Y \\
\searrow && \swarrow  \\
&Z& \\
\end{array}$$
such that
\begin{enumerate}
\item $X = {\Bbb C}^x \times Z$;
\item $Y = {\Bbb C}^y \times Z$;
\item $x \ge y$;
\item $Y$ is rational.
\end{enumerate}
Then $X$ is rational.
\end{lem}

\proof $X = {\Bbb C}^x \times Z \cong {\Bbb C}^{x - y}  \times 
({\Bbb C}^y \times Z) = {\Bbb C}^{x - y}  \times Y$.
\endproof

\begin{say} Consider a diagram
$$\begin{array}{ccccccccccccc}
     (n;d) &&(n_1;d_1) &&& \ldots && (n_{t-1}; d_{t-1}) &&(n_t; d_t) \\
       \searrow && \swarrow  \searrow &&& \ldots &&  \swarrow \searrow &&  \swarrow \\
            &(\ell_1; k_1)&& (\ell_2; k_2) && \ldots &&  & (\ell_t; k_t) 
\end{array} \eqno (6.8) $$
where each pair in the diagram is a pair of positive integers satisfying (6.1)
and each downward arrow represents  successive reductions and dual reductions.
For example, $(\ell_1; k_1)$ is obtained from $(n; d)$ and $(n_1; d_1)$
 by some sequences of reductions and dual reductions, respectively. 
It is allowed that $(n_t; d_t)=(\ell_t; k_t)$.

Correspondingly, we get a diagram of dominant rational maps:
$$\begin{array}{ccccccccccccc}
     {\cal M}(n,L)&&{\cal M}(n_1, L_1)&& \ldots && {\cal M}(n_t, L_t) \\
          \searrow && \swarrow  \searrow && \ldots &&   \swarrow \\
             &{\cal M}(\ell_1, K_1)&& \cM ( (\ell_2; K_2) & \ldots  & {\cal M}(\ell_t, K_t)
\end{array} $$
where $L_i$ is a line bundle with deg$L_i=d_i$ and $K_i$ is a line bundle
with deg$K_i=k_i$. 
\end{say}

\begin{defn}
\label{admissible}
The diagram (6.8) is called admissible if 
$n \ge n_i$ and $\gcd (l_i, k_i) = 1$ for all $i$.
\end{defn}

\begin{thm} 
\label{ballicothm}
Suppose the diagram (6.8) is admissible. Then 
${\cal M}(n, L)$ is birational to ${\Bbb C}^{z_t}\times {\cal M}(n_t,L_t)$
where $z_t=(n^2-n^2_t)(g-1)\ge 0$ is the difference of the dimension of
two moduli spaces.
\end{thm}

\proof We shall use Proposition
\ref{rationalityinductioncor} and Lemma \ref{transferringlemma} repeatedly.
Let $x= \dim {\cal M}(n, L)$, $x_i = \dim {\cal M}(n_i, L_i)$, and
$y_i = \dim {\cal M}(l_i, K_i)$. For simplicity, we use the notation
$X \cong^{bir} Y$ to indicate that $X$ and $Y$ are birational equivalence.
Then we have
$${\cal M}(n, L) \cong^{bir} {\Bbb C}^{x - y_1} \times \cM (l_1, K_1)
= {\Bbb C}^{x - x_1} \times {\Bbb C}^{x_1 - y_1} \times \cM (l_1, K_1) $$
$$ \cong^{bir} {\Bbb C}^{x - x_1} \times \cM (n_1, L_1)  \cong^{bir} 
 {\Bbb C}^{x - x_1} \times{\Bbb C}^{x_1 - y_2} \times \cM (l_2, K_2) $$
$$ = {\Bbb C}^{x - y_2} \times \cM (l_2, K_2) =  
{\Bbb C}^{x - x_2} \times{\Bbb C}^{x_2 - y_2} \times \cM (l_2, K_2) $$
$$\cong^{bir} {\Bbb C}^{x - x_2} \times \cM (n_2, L_2) \cong^{bir}
 \ldots  \cong^{bir} {\Bbb C}^{x - x_t} \times \cM (n_t, L_t).$$
The theorem then follows.
\endproof 

\begin{defn}
\label{ballicopair}
A pair $(n; d)$ is called a fine pair if there exists an admissible diagram (6.8)
such that $(n_t; d_t)$ is a nice pair. Fine pairs are also necessarily coprime.
\end{defn}

\begin{thm}
\label{ballicopairthem}
If $(n; d)$ is a fine pair, then ${\cal M}(n, L)$ is rational.
\end{thm}

\proof It follows from Theorem \ref{ballicothm} and Theorem \ref{newsteadpairthm}.
\endproof

\begin{exmp}
Assume $g=6$. Then the  pair $(7; 38)$ is not a nice pair by Corollary \ref{iterationcor}
($\gcd(38, 6)=2$ and $\gcd(38+7, 6)=3$). Notice that $\gcd(7, 38)=1$.

One can check the following statements:
\begin{enumerate}
\item $(7; 38)$ is the (one-step) dual reduction of $(11+7m; 62+35m)$ for all integer $m\ge 0$.
\item  $m=7$ is the smallest number such that the pair $(11+7m; 62+35m)$ is a nice pair. 
        (When $m=7$, the pair is $(60; 307)$ and we actually have
        $$(60; 307)\longrightarrow (53; 307)\longrightarrow (11; 65) \longrightarrow (1; 6)$$
       where all $``\longrightarrow "$ are (one-step) reductions.)
\item  When $m$ is even, $(11+7m; 62+35m)$ is not a nice pair.
       (The dual reduction of the pair gives $(7;38)$ which is not a nice pair; a reduction gives
          $(4+7m; 62+35m)$ which  is not a nice pair either since
           $2$ divides $\gcd(4+7m, 62+35m)$ when $m$ is even.)  
\end{enumerate}

However $(11+7m; 62+35m)$ are fine pairs for all $m\ge 7$:
$$\begin{array}{clc}
(11+7m; 62+35m) && (60; 307)  \\
\searrow && \swarrow \\
&(7; 38)& \\
\end{array}$$
Hence  the corresponding
moduli spaces ${\cal M}(11+7m, L)$ where ${\rm deg}L=62+35m$ are rational
for all $m\ge 7$ by Theorem \ref{ballicopairthem}.

The above provides infinitely many fine pairs that are not nice pairs. 
\end{exmp}

\begin{rem} Finally,  some remarks are in order.
\begin{enumerate}
\item There are coprime pairs that are not fine pairs.
\item In [{\bf N2}], Newstead defined a good pair to be
a coprime pair $(n; d)$ whose corresponding moduli space is
rational. Nice pairs and fine pairs are good. The converse
may not be true. To find all good pairs (that are not nice or fine)
seems requiring methods different than the one explored in this paper.
\item There is an algorithm to locate nice and fine pairs on the lattice cone in the $xy$-plane
       defined by inequalities:
      $$ x(g-1) < y \le xg$$ for any fixed $g \ge 2$. A reduction takes the form
$$(x, y) \rightarrow (gx -y, y - k (gx-y))$$
   for some non-negative integer $k$. A dual reduction takes the form
$$(x, y) \rightarrow (y - (g-1) x, (2g-1)x -y - k (y-(g-1)x))$$ for some non-negative integer $k$.
\item By using the techniques of variations of moduli spaces of parabolic
bundles, H. Boden and K. Yokogawa were able to prove the rationality of 
certain moduli spaces. (see [{\bf B, BH}].
\end{enumerate}
\end{rem}

%\bibliographystyle{amsplain}
%\makeatletter \renewcommand{\@biblabel}[1]{\hfill#1.}\makeatother
%\newcommand{\bysame}{\leavevmode\hbox to3em{\hrulefill}\,}
%\begin{thebibliography}{10}
\noindent
{\bf Reference}

\vskip 1cm

%\begin{enumerate}
\noindent
{[{\bf BPS}]} C. Banica, M. Putinar, G. Schumacher,
Variation der globalen Ext in Deformationen kompakter
komplexer R\"aume
{\it Math. Ann.}
 {\bf 250} (1980), 135-155

\noindent
{[{\bf Ba}]} E. Ballico,
Stable rationality for the variety of vector bundles over
an algebraic curve,
{\it J. London Math. Soc.} {\it 2nd Series} {\bf 30} (1984), 21-26.

\noindent [{\bf B}]
H. Boden,
Rationality of the moduli space of vector bundles over a smooth curve,
preprint IHES/M/95/73.

\noindent [{\bf BY}] H. Boden, K. Yokogawa,
Rationality of moduli spaces of parabolic bundles,
alg-geom/9610013

\noindent
{[{\bf B1}]} J.E. Brosius,
 Rank-2 vector bundles on ruled surfaces, I,
{\it Math. Ann.}
{\bf 265} (1983), 155-168.

\noindent
{[{\bf B2}]} J.E. Brosius,
 Rank-2 vector bundles on ruled surfaces, II,
{\it Math. Ann.}
{\bf 266} (1983), 199-214.

%\noindent   {[{\bf D}]} I. Dolgachev,
%Rationality of fields of invariants,
%in {\it Algebraic Geometry}, Bowdoin (1985), AMS (1987).

%\noindent   {[{\bf DH}]} I. Dolgachev and Y. Hu,
%Variation of Geometric Invariant Theory Quotients,
%{\it Publ. Math. I.H.E.S.} (to appear).

%\noindent   {[{\bf DO}]} I. Dolgachev and Ortland,
%Point sets in projective space and Theta functions,
%{\it Ast\'erisque, Soc. Math. de France}
% {\bf  68} (1992).

%\noindent   {[{\bf F}]} R. Friedman, The geometry of vector bundles over algebraic varieties, (to appear).

\noindent   {[{\bf GH}]} P. Griffiths and J. Harris,
Principles of algebraic geometry,
Wiley, New York (1978).

\noindent   {[{\bf G}]} I. Grzegorczyk, On Newstead's conjecture on vector bundles
on algebraic curves, {\it Math. Ann.}  {\bf 300} (1994), 521-541.

%\noindent   {[{\bf Ka}]} P. Katsylo, Rationality of the orbit spaces of irreducible representations
%of the group $SL_2$,
%{\it Izv. Akad. Nauk SSSR Ser. Mat.}
% {\bf  47} (1983), 26-36.

\noindent   {[{\bf L}]} H. Lange,
Universal families of extensions,
{\it  Journal of Algebra}
{\bf  83} (1983), 101-112.

\noindent   {[{\bf M}]} M. Maruyama, Elementary transformationations in the
theory of algebraic vector bundles, {\it LNM} {\bf 961} (1982), 241-266.

%\noindent   {[{\bf MS}]} V. Mehta and C. Seshadri,
% Moduli of vector bundles with parabolic structures over curves,
%{\it   Math. Ann.}, {\bf 101,} (1980), 205-239.

\noindent   {[{\bf MF}]} D. Mumford and J. Fogarty,
   Geometric Invariant Theory,
{\it A Series of Modern Surveys in Mathematics}, Springer-Verlag (1982).

\noindent   {[{\bf NR}]} M. S. Narasimhan and S. Ramanan: Moduli of vector bundles
on a compact Riemann surface, {\it Ann. of Math.} {\bf 82} (1965), 540-567.

\noindent   {[{\bf N1}]} P. Newstead,
Rationality of moduli spaces of stable bundles,
{\it Math. Ann. } {\bf 215} (1975),
251-268. 

\noindent   {[{\bf N2}]} P. Newstead,
Rationality of moduli spaces of stable bundles,
{\it Correction,} {\it Math. Ann. } {\bf 249} (1980),
281-282.

\noindent   {[{\bf R}]} S. Ramanan: The moduli spaces of vector bundles
over an algebraic curve, {\it Math. Ann. } {\bf 200} (1973), 69-84.
%\noindent   {[{\bf Sh}]} N. Shepherd-Barron,
%Rationality of moduli spaces via invariant theory,
%in {\it Topological methods in algebraic transformation  groups},
%Birkhauser, 1989. 

%\noindent {[{\bf NS}]} M. S. Narasimhan  and C. S. Seshadri, Stable
%and unitary vector bundles on a compact Riemann surface,  {\it Ann. Math.} 
%{\bf 82} (1965) 540-567.

\noindent   {[{\bf T1}]} A. Tjurin,
Classification of vector bundles over an algebraic curve
of arbitrary genus,
{\it Amer. Math. Soc. Transl.} {\bf 63}
(1967),  245-279.

\noindent   {[{\bf T2}]} A. Tjurin,
Classification of $n$-dimensional vector bundles over an algebraic curve
of arbitrary genus,
{\it Amer. Math. Soc. Transl.} {\bf 73}
(1968),  196-211.

\noindent   {[{\bf T3}]} A. Tjurin, On the classification of two-dimensional
vector bundles over an algebraic curve of arbitrary genus, {\it
Izv. AKa. Nauk SSSR Ser. Mat.} {\bf 28} (1964), 21-52.

%\noindent   {[{\bf Se}]} J. P. Serre, 
%Les espaces fiber de principal,
%{Seminaire C. Chevalley} vol 2, 1958. 

%\end{enumerate}

%\end{thebibliography}

\vskip .5cm
\noindent

\noindent
Y.H.   Department of Mathematics, University of California, Berkeley, 
CA 94720, USA. hu@@math.berkeley.edu

%\vskip .4cm
%\noindent
%Department of Mathematics, University of Texas, Arlington, TX 76019.

\vskip .5cm

\noindent
W.L.   Department of Mathematics, Hong Kong University of Science and Technology, 
     Clear Water Bay, HK.
    mawpli@@uxmail.ust.hk

\end{document}